\documentclass[final,3p,times]{elsarticle}

\usepackage{amssymb}
\usepackage{amsthm}
\usepackage{amsmath} 

\usepackage{lineno}

\journal{Journal of Sound and Vibration}

\begin{document}

\begin{frontmatter}




\title{A stable decoupled perfectly matched layer for the 3D wave equation using the nodal discontinuous Galerkin method}

\author[label1]{Sophia Julia Feriani}
\author[label2]{Matthias Cosnefroy}
\author[label3]{Allan Peter Engsig-Karup}
\author[label4]{Tim Warburton}
\author[label2]{Finnur Pind}
\author[label1]{Cheol-Ho Jeong}

\affiliation[label1]{organization={Acoustic Technology, Department of Electrical and Photonics Engineering, Technical University of Denmark},
            addressline={Ørsteds Plads, Bygning 352}, 
            city={Kongens Lyngby},
            postcode={2800}, 
            state={},
            country={Denmark}}
\affiliation[label2]{organization={Treble Technologies},
            addressline={Hafnartorg, Kalkofnsvegur 2}, 
            city={Reykjavík},
            postcode={101}, 
            state={},
            country={Iceland}}      
\affiliation[label3]{organization={Scientific Computing, Department of Applied Mathematics and Computer Science, Technical University of Denmark},
            addressline={Richard Petersens Plads, Building 324}, 
            city={Kongens Lyngby},
            postcode={2800}, 
            state={},
            country={Denmark}}
\affiliation[label4]{organization={Data \& Decision Sciences, Virginia Tech
},
            addressline={727 Prices Fork Road}, 
            city={Blacksburg},
            postcode={24060}, 
            state={VA},
            country={USA}}
    
\begin{abstract}
In outdoor acoustics, the calculations of sound propagating in air can be computationally heavy if the domain is chosen large enough to fulfil the Sommerfeld radiation condition. By strategically truncating the computational domain with a efficient boundary treatment, the computational cost is lowered. One commonly used boundary treatment is the perfectly matched layer (PML) that dampens outgoing waves without polluting the computed solution in the inner domain.
The purpose of this study is to propose and assess a new perfectly matched layer formulation for the 3D acoustic wave equation, using the nodal discontinuous Galerkin finite element method. The formulation is based on an efficient PML formulation that can be decoupled to further increase the computational efficiency and guarantee stability without sacrificing accuracy. This decoupled PML formulation is demonstrated to be long-time stable and an optimization procedure of the damping functions is proposed to enhance the performance of the formulation.
\end{abstract}

\begin{keyword}
outdoor sound propagation \sep nodal discontinuous Galerkin method \sep perfectly matched layer \sep time domain
\end{keyword}

\end{frontmatter}


\section{Introduction} \label{sec:intro}
Outdoor acoustics deals with sound propagation in open environments and is of high practical relevance in areas such as urban planning, noise control, environmental protection, architectural design, event planning, and public health. When simulating outdoor sound propagation, the domain of interest is typically confined spatially and the computation resources are limited. Therefore, to use the computational resources efficiently, it is often done by limiting the size of the computational domain. However, the truncation of an otherwise open domain may introduce unwanted so-called spurious reflections that can unintentionally pollute the numerical solution if the physical properties of the system are not modelled carefully.

To mitigate unintended effects of domain truncation, it is crucial to design and employ non-reflecting boundary treatments. Non-reflective boundary conditions \cite{Bayliss1980,Engquist1979,Collino1993} and non-reflective boundary layer \cite{Appelo2009,Israeli1981,Shlomo1995} approaches are typically used. The former approach imposes a condition on the outer domain boundary to minimize the spurious reflections coming back into the domain of interest, however, they are not simple to implement if high accuracy is required. The latter approach alter the governing equations in a buffer zone, typically wide, to ensure a minimum reflection at the cost of large additional computational demand. The perfectly matched layer boundary treatment is commonly used due to its ease of implementation, accuracy, and versatility and the main goal of this study is to develop a stable perfectly matched layer for the acoustic wave equation using the nodal discontinuous Galerkin method.

The Perfectly Matched Layer (PML) technique introduces a change of material at the periphery of the outer boundary of the domain where the material properties are designed to damp propagating waves \cite{Berenger1994}. PMLs can be interpreted as a complex coordinate stretching of the governing equations in the frequency domain \cite{Chew1994} and it is designed to have no reflection for any waves at the continuous level. However, once the problem is discretized numerically, some issues arise in terms of well-posedness \cite{Abarbanel1997}, stability \cite{Hesthaven1998,Becache2003} and accuracy \cite{Diaz2006} (especially for grazing incidence waves) which prompted ongoing research and the creation of different formulations of PML in the last decades \cite{Hu2001,Hu2005,Parrish2009,Meza2008,Roden2000,Grote2010,Etienne2010}.

In the field of acoustics, Hu and colleagues proposed a ``uni-axial'' PML \cite{Hu2001,Hu2005,Parrish2009} 
, however, this formulation can generate late-time instabilities. Those instabilities can be fixed if a stabilising correction terms in the complex coordinate stretching is included. This formulation called "multi-axial" PML \cite{Meza2008} generalize the U-PML by adding an orthogonal term to the damping function but the layer is no longer perfectly matching. In theory, PML is "perfectly matched" to the governing equations meaning that waves reaching the PML from all angles should be absorbed. After discretization, grazing incidence waves are not fully absorbed and can generate spurious reflections spreading to the inner domain. The "convolutional" PML (C-PML) was formulated to address this issue \cite{Komatitsch2007,Martin2009} and the convolution can be made more computationally efficient by calculating it in a recursive way. However, this PML formulation introduces a high number of auxiliary variables and corresponding PDEs to be solved in parallel to the main governing equations. Moreover, the C-PML was shown to be not that efficient at absorbing grazing incidence waves \cite{Cosnefroy2019}. 

Later, an efficient PML formulation for the 3D wave equation in its standard second-order form was proposed \cite{Grote2010} that involves fewer auxiliary variables and corresponding PDEs than its counterparts. Moreover, stability was proven for a PML in a single direction \cite{Grote2010,Becache2021,Baffet2019} for a constant damping but not in all three directions. This formulation was implemented for the second-order wave equation using the Finite Difference Time Domain method (FDTD) \cite{Grote2010,Becache2021} and Finite Element Method (FEM) \cite{Kaltenbacher2013}. 

In the recent two decades, the Discontinuous Galerkin Finite Element Method (DG-FEM) \cite{Hesthaven2007} has seen steadily growing interest in advanced applications, and in particular for the solution of hyperbolic partial differential equations utilising complex and adaptive mesh representations for the numerical discretization. The numerical method was first proposed for the neutron transport equation in the 1972 \cite{ReedHill1972} and is considered particularly useful in room and outdoor acoustics \cite{Atkins1998,Toulopoulos2006} because it is high-order accurate that is beneficial for long-term wave propagation problems \cite{Kreiss1972} and in the field of acoustics, realistic problems involve complex geometries requiring unstructured grid that can be curved or adapted with local mesh refinement \cite{Langtangen1998}. Moreover, unlike FDTD which is a method usually implemented for outdoor acoustics, DG-FEM is highly parallelizable \cite{Melander2023} since the governing equations can be integrated in an element-wise manner only requiring solution traces from elements sharing a face \cite{Hesthaven2007}.

In the field of acoustics, the nodal DG-FEM has been used to simulate sound behavior in rooms \cite{Wang2019}, including for different types of reflective boundary condition \cite{Wang2020,Pind2020,Pind2021}. 
In outdoor acoustics, \cite{Modave2017} constructed two PML formulation for convex domains with regular curved boundary and used the discontinuous Galerkin method to discretize the problem. Those formulations are however not simple to implement and are limited to specific mesh configurations. 
In the field of geophysics, \cite{Etienne2010} designed a convolutional PML using DG-FEM for the first-order elastic wave equation but instabilities appear for an isotropic medium. Those instabilities can be mitigated by using a wider PML absorbing region with a corresponding increase in the computation time.  

This paper proposes a stable and efficient PML for dealing with open boundaries when solving the 3D acoustic wave equations discretized with the nodal DG-FEM. Since the formulation of the PML of \cite{Grote2010} is stable for a damping from a single direction, the idea in this study is to rely on three stable and independent PML formulations for each Cartesian direction and combined them to form what we refer to as a {\em decoupled PML formulation}. The resulting formulation involves fewer auxiliary variables and corresponding PDEs than other PML formulations, which make it computationally efficient to solve. The damping performance, accuracy and stability of the PML formulation is assessed.

The organisation of this paper is as follows. Section \ref{sec:NumericalMethod} presents the numerical method used to solve the governing equations and the modifications needed to incorporate the PML formulation. Section \ref{sec:method} describe the method for the results of section \ref{sec:Results}. Finally, those results are discussed in section \ref{sec:Discussion}.

\section{Numerical method} \label{sec:NumericalMethod}
\subsection{Governing equations in air} \label{sec:GovEq}
The governing equations for acoustic wave propagation in a homogeneous motionless medium take the form of the first-order partial differential equations as follows:
\begin{subequations} \label{1storder}
\begin{align}
    \frac{\partial p}{\partial t} &+ \rho c^2 \nabla \cdot \textbf v = 0, \label{1storderP}\\
    \frac{\partial \textbf v}{\partial t} &+ \frac{1}{\rho} \nabla p = \textbf 0, \label{1storderV}
\end{align}
\end{subequations}
where $p = p(\textbf x,t)$ is the acoustic pressure [$Pa$] and $\textbf v =\textbf v(\textbf x,t)$ is the particle velocity [$m/s$] with $\textbf x = [x,y,z]^T$ the position in space in the domain $\Omega \subset \mathbb{R}^3$ and $t$ the time [$s$]. The speed of sound $c$ is set to $343 \; m/s$ and the air density $\rho$ is set to $1.2 \; kg/m^{3}$ in this study. Eq. \eqref{1storder} can be written as the following hyperbolic system:
\begin{subequations}
\begin{equation} \label{1storderMatrix}
    \frac{\partial \textbf u}{\partial t} = - A_x \frac{\partial \textbf u}{\partial x} - A_y \frac{\partial \textbf u}{\partial y} - A_z \frac{\partial \textbf u}{\partial z},
\end{equation}
where 
\begin{equation}
A_x = \begin{bmatrix} 0& \rho c^2 &0&0 \\ 1/\rho&0&0&0 \\ 0&0&0&0 \\ 0&0&0&0 \end{bmatrix}, \quad A_y = \begin{bmatrix} 0 & 0 & \rho c^2 & 0 \\ 0&0&0&0 \\ 1/\rho&0&0&0 \\ 0&0&0&0 \end{bmatrix}, \quad A_z = \begin{bmatrix} 0 & 0 & 0 & \rho c^2 \\ 0&0&0&0 \\ 0&0&0&0 \\ 1/\rho&0&0&0 \end{bmatrix}.
\end{equation}
\end{subequations}
and $\textbf{u} = \textbf u(\textbf x ,t) = [p,v_x,v_y,v_z]^T$ is the matrix of the acoustic variables at time $t$.

\subsection{Spatial discretization}\label{sec:SpatialDiscrete}
The first-order PDEs Eq. \eqref{1storder} are spatially discretized using DG-FEM following \cite{Hesthaven2007}.
The three-dimensional spatial domain $\Omega$ is approximated by a computational domain $\Omega_h$, partitioned into $K$ non-overlapping tetrahedral elements of volume $V^k$ so $\Omega \approx \Omega_h = \cup^K_{k=1}V^k$. The global solution $\textbf u(\textbf x ,t)$ is approximated by the piece-wise $N$-th order polynomial approximation $\textbf u_h(\textbf x,t)$, which is the direct sum of the $K$ local polynomial solutions $\textbf u_h^k(\textbf x,t)$ with $k=\{1,...,K\}$. In the element $k$ of volume $V^k$ in the tessellation of the domain, the local approximate solution $\textbf u^k_h(\textbf x,t)$ is:
\begin{equation}
    \textbf u^k(\textbf x,t) \approx \textbf u^k_h(\textbf x,t) = \sum^{N_p}_{i=1} \textbf u^k_h(\textbf x_i,t) \varphi_i^k(\textbf x),
\end{equation}
where $\textbf u^k_h(\textbf x_i,t)$ are the acoustic variables at the node $\textbf x_i$ and time $t$ in element $k$.
$\varphi_i^k$ are the nodal basis functions and are taken as the interpolating nodal Lagrangian polynomial basis function associated to the $i$-th node at position $\textbf x_i$ in each element $k \in \{1,...,K\}$.
$N_p$ is the number of degrees of freedom per element and is equal to $N_p = (N+3)(N+2)(N+1)/6$ in 3D. 
In the following, the subscript $h$ referring to the numerical approximation to the solution in the mesh is dropped to make the notation simpler.

Starting from Eq. \eqref{1storderMatrix}, it is multiplied by this basis function and integrated over the volume of element $k$. After a first integration by parts, the weak formulation can be expressed as:
\begin{equation}
    \int_{V^k} \varphi_i^k \frac{\partial \textbf u^k}{\partial t}  d \textbf x= \int_{V^k} \bigg( \frac{\partial \varphi_i^k}{\partial x}  A_x \textbf u^k + \frac{\partial \varphi_i^k}{\partial y} A_y \textbf u^k + \frac{\partial \varphi_i^k}{\partial z}  A_z \textbf u^k  \bigg)  d \textbf x - \int_{\partial V^k}  \varphi_i^k \; \textbf n  \cdot \textbf F^*(\textbf u^k, \textbf u^{k,ext}) dS,
\end{equation}
where $\textbf n = [n_x,n_y,n_z]^T$ is the local normal vector pointing outward at the boundary element and $\textbf F^*$ is the numerical flux between the interior element '$k$' and the neighbouring element identified with the subscript '${k,ext}$', introduced because the solution is discontinuous at the interfaces between adjacent elements. After a second integration by parts, the strong formulation is expressed as:
\begin{equation}\label{StrongForm}
\int_{V^k} \varphi_i^k \frac{\partial \textbf u^k}{\partial t}  d\textbf x= 
- \int_{V^k} \varphi_i^k \bigg(  A_x \frac{\partial \textbf u^k }{\partial x} + A_y \frac{\partial \textbf u^k}{\partial y} + A_z \frac{\partial \textbf u^k }{\partial z} \bigg)  d \textbf x
+ \int_{\partial V^k} \varphi_i^k (A_n \textbf u^k - \textbf n \cdot \textbf F^*(\textbf u^k, \textbf u^{k,ext})) dS,
\end{equation}
where $A_n = A_x n_x + A_y n_y + A_z n_z$.
To compute the volume integral, the local mass matrix $\mathcal{M}_{ij}^k$ and local stiffness matrix $(\mathcal{S}^k)_{ij}$ are defined for the element $k$ as:
\begin{equation}
    \mathcal{M}_{ij}^k = \int_{V^k} \varphi_i^k(\textbf x) \varphi_j^k(\textbf x) d \textbf x, \quad (\mathcal{S}^k)_{ij} = \int_{V^k} \varphi_i^k(\textbf x) \frac{d \varphi_j^k(\textbf x)}{d\theta} d\textbf x.
\end{equation}
To compute the surface integral, $\mathcal E$ is introduced. It is the matrix operator along the four faces of the tetrahedron element $k$, so the strong form Eq. \eqref{StrongForm} becomes:
\begin{equation} \label{NewStrongForm}
    \mathcal{M}^k \frac{\partial \textbf u^k}{\partial t} = - A_x \mathcal S_x^k  \textbf u^k - A_y \mathcal S_y^k \textbf u^k - A_z \mathcal S_z^k \textbf u^k + \mathcal{E} \;  \textbf n \cdot (A_n \textbf u^k -\textbf F^*(\textbf u^k, \textbf u^{k,ext})).
\end{equation}

An upwind numerical flux is taken to calculate the surface integral in the last term of Eq. \eqref{NewStrongForm}. The eigen-decomposition of $A_n$ provides the four eigenvalues stored in the diagonal of a matrix $\Lambda$ and the corresponding right eigenvectors are stored in the columns of a matrix $W$. Since Eq. \eqref{1storderMatrix} is a constant-coefficient linear problem, Roe's method can be applied \cite{LeVeque2002} so the upwind numerical fluxes at the discontinuity between all elements are written as:
\begin{equation} \label{NumFlux}
    \textbf F^*(\textbf u^k, \textbf u^{k,ext}) = \frac{1}{2} \left(  A_n \textbf u^k + A_n \textbf u^{k,ext} + W|\Lambda|W^{-1}(\textbf u^k - \textbf u^{k,ext})\right).
\end{equation}
At the termination of the mesh, one or more faces of the border tetrahedral elements have no neighbours. Two types of boundary condition are therefore imposed on $\textbf u^{k,ext}$ to compute the numerical flux Eq. \eqref{NumFlux}: either perfectly reflective ($p^{k,ext} = p^k$ and $\textbf v^{k,ext} = -\textbf v^k$) or a simple absorbing boundary condition ($\textbf u^{k,ext} = \textbf 0$) called "ABC" in the rest of the paper. 

\subsection{Temporal discretization}\label{sec:TemporalDiscrete}
An explicit low-storage forth-order Runge-Kutta method is used to integrate the semi-discrete coupled equations in time \cite{Carpenter1994}. The PDE Eq. \eqref{1storderMatrix} is the semi-discrete expression that can be re-written as a system of ODEs in a general form:
\begin{equation}
    \frac{d \textbf u_h}{d t}=\mathcal{L}_h(\textbf u_h(t),t),
\end{equation}
where $\textbf u_h$ is the vector of unknown variables and $\mathcal{L}_h$ a spatial operator which corresponds to the right-hand side of Eq. \eqref{1storderMatrix}. The five stages of the scheme are calculated as follows:
\begin{subequations}
\begin{align}
&\textbf u^{(0)}_h = \textbf u^{n}_h, \\
&\alpha\in\{1,...,5\} \; \begin{cases}
    k^{(\alpha)} = a_\alpha k^{(\alpha-1)}+\Delta t \mathcal{L}_h(\textbf u_h^{(\alpha-1)},t_n), \\
    \textbf u^{(\alpha)}_h = \textbf u_h^{(\alpha-1)} + b_\alpha k^{(\alpha)},
\end{cases} \\
&\textbf u_h^{n+1}=\textbf u_h^{(5)},
\end{align}
\end{subequations}
where $\Delta t = t^{n+1}-t^n$ is the time step [$s$] and the value of the coefficients $a_\alpha$ and $b_\alpha$ are found in \cite{Carpenter1994}. The time step $\Delta t$ depends on the size of the mesh elements and the order of approximation $N$ and is calculated as:
\begin{equation} \label{timeStepDef}
    \Delta t = CFL \frac{\min_{k}(\Delta x^k)}{c}\frac{1}{N^2},
\end{equation}
where $CFL=\mathcal{O}(1)$ is the Courant–Friedrichs–Levy number taken as $0.95$ and  $\min_{k}(\Delta x^k)$ is the minimal characteristic length among all the elements of the mesh. The $1/N^2$ scaling stems from the use of node-optimized non-uniform distributions \cite{Warburton2008} across the elements in the mesh.

\subsection{Perfectly Matched Layer}
In this subsection, the formulation of the PML is constructed in a region surrounding the 3D domain of interest $\Omega$ called $\Omega^{pml}$ containing $K^{pml}$ elements. The width of the PML region for the $x$-$y$-$z$-directions is $\delta_x$, $\delta_y$ and $\delta_z$ respectively. 
The continuous PML formulation is then discretized to fit into the nodal DG-FEM framework in the time domain described in the previous subsections \ref{sec:SpatialDiscrete} and \ref{sec:TemporalDiscrete}.
\subsubsection{Continuous formulation}
The first-order PDE system Eq. \eqref{1storder} can be reduced to a single second-order equation: 
\begin{equation}
    \frac{\partial^2p}{\partial t ^2} = p_{tt} = c^2 \, \nabla\cdot\nabla p.
\end{equation}
In this section, the subscript '$t$' correspond to the first derivative in time and '$tt$' the second derivative in time. Let's recall the Laplace transform of $p$ from the time domain to the frequency domain:
\begin{equation}
    \hat p = \hat p(\textbf x,s)=\int^\infty_0 e^{st}p(\textbf x,t) dt,
\end{equation}
where $s \in \mathbb{C}$ is the complex frequency variable [$s^{-1}$]. At initial time, the acoustic variables are chosen so that they are null in the PML region ; so the Laplace transform $\hat p$ satisfies: 
\begin{equation} \label{LaplaceOfHelm}
    s^2\hat p = c^2 \frac{\partial}{\partial x} \left( \frac{\partial \hat p}{\partial x} \right) + c^2\frac{\partial}{\partial y} \left( c^2 \frac{\partial \hat p}{\partial y} \right)+ c^2\frac{\partial}{\partial z} \left( \frac{\partial \hat p}{\partial z} \right).
\end{equation}

The role of the PML is to dampen the wave propagating inside $\Omega^{pml}$ therefore three continuous, smooth and constant in time damping functions are defined along each Cartesian direction, they are denoted $\sigma_x(x)$, $\sigma_y(y)$ and $\sigma_z(z)$. As in all other PML formulations previously cited, the simplest complex coordinate stretching transformation \cite{Chew1994} for $x_i = \{ x,y,z \}$ is defined as:
\begin{equation} \label{coordStretch}
    \frac{\partial}{\partial x_i} \to \frac{\partial}{\partial \tilde x_i} = \frac{s}{s+\sigma_i} \frac{\partial}{\partial x_i} = \frac{1}{1+\frac{\sigma_i}{s}} \frac{\partial }{\partial x_i} = \frac{1}{\gamma_i} \frac{\partial }{\partial x_i},
\end{equation}
where $\gamma_i=\gamma_i(\sigma_i,s)$ scales the spatial partial derivative in the frequency domain. In $\Omega^{pml}$, $\hat p$ satisfies: 
\begin{equation}
    s^2\hat p =c^2 \frac{\partial}{\partial \tilde x} \left( \frac{\partial \hat p}{\partial \tilde x} \right) + c^2\frac{\partial}{\partial \tilde y} \left( \frac{\partial \hat p}{\partial \tilde y} \right)+ c^2\frac{\partial}{\partial \tilde z} \left( \frac{\partial \hat p}{\partial \tilde z} \right).
\end{equation}
It has been proven that the PML formulation of \cite{Grote2010} along only one of the three Cartesian directions is stable for a constant damping profile \cite{Grote2010,Becache2021,Baffet2019}. Therefore, let's first consider the case where only one damping function is not null, for example, $\sigma_x \neq 0$ and $\sigma_y = \sigma_z = 0$:
\begin{subequations}
\begin{equation}
    \gamma_x s^2\hat p = \frac{1}{\gamma_x} c^2 \frac{\partial^2 \hat p}{\partial x^2} + \gamma_x c^2 \frac{\partial^2 \hat p}{\partial y^2} + \gamma_x c^2 \frac{\partial^2 \hat p}{\partial z^2}.
\end{equation}
By replacing $\gamma_x$ with its definition:
\begin{equation}
    \left(1+\frac{\sigma_x}{s}\right) s^2\hat p = \left(1-\frac{\sigma_x}{s+\sigma_x}  \right) c^2 \frac{\partial^2 \hat p}{\partial x^2} + \left(1+\frac{\sigma_x}{s}\right) c^2 \frac{\partial^2 \hat p}{\partial y^2} + \left(1+\frac{\sigma_x}{s}\right) c^2 \frac{\partial^2 \hat p}{\partial z^2},
\end{equation}
and after re-arranging the term, we get:
\begin{equation}
    s^2\hat p + \sigma_x s\hat p = c^2 \Delta p
    -\frac{\sigma_x}{s+\sigma_x} c^2 \frac{\partial^2 \hat p}{\partial x^2} + 
    \frac{\sigma_x}{s} c^2\frac{\partial^2 \hat p}{\partial y^2} + \frac{\sigma_x}{s} c^2 \frac{\partial^2 \hat p}{\partial z^2} , \label{sPres}
\end{equation}
\end{subequations}
Following \cite{Grote2010}, three unknown variables $\boldsymbol{\phi} = \boldsymbol{\phi}(\textbf x,t) = [\phi^x,\phi^y,\phi^z]^T$ are introduced to come back to the time domain formulation. Those auxiliary variables must solve:
\begin{equation}
    \nabla \cdot \boldsymbol{\phi} = -\frac{\sigma_x}{s+\sigma_x} c^2 \frac{\partial^2 \hat p}{\partial x^2} + \frac{\sigma_x}{s} c^2\frac{\partial^2 \hat p}{\partial y^2} + \frac{\sigma_x}{s} c^2 \frac{\partial^2 \hat p}{\partial z^2},
\end{equation}
which can be re-arranged as:
\begin{subequations} \label{sphi}
\begin{align}
    \phi^x &= -\frac{\sigma_x}{s+\sigma_x} c^2 \frac{\partial \hat p}{\partial x}, \\
    \phi^y &= \frac{\sigma_x}{s} c^2\frac{\partial \hat p}{\partial y}, \\
    \phi^z &= \frac{\sigma_x}{s} c^2\frac{\partial \hat p}{\partial z},
\end{align}
\end{subequations}
Equations \eqref{sPres} and Eq. \eqref{sphi} are brought back to the time domain so the PML formulation with damping only in the $x$-direction is:
\begin{subequations} \label{PML2norder}
\begin{align}
    p_{tt}+\sigma_x p_t &= c^2 \Delta p +\nabla \cdot  \boldsymbol{\phi}, \label{PML2nordera}\\
    \boldsymbol{\phi} _t &= \Gamma_1^x  \boldsymbol{\phi} + c^2 \Gamma_2^x \nabla p, \label{PML2norderb}
\end{align}
where 
\begin{equation}
\Gamma_1^x=\begin{bmatrix}-\sigma_x&0&0\\0&0&0\\0&0&0\end{bmatrix}, \quad  \Gamma_2^x=\begin{bmatrix}-\sigma_x&0&0\\0&\sigma_x&0\\0&0&\sigma_x\end{bmatrix}.
\end{equation}
\end{subequations}
To discretize Eq. \eqref{PML2nordera} in the same framework as given in Section \ref{sec:SpatialDiscrete}, it first needs re-writing into a first order form. By taking the temporal derivative of the particle velocity as $\textbf v_t = -\frac{1}{\rho} \nabla p - \frac{1}{\rho c^2} \boldsymbol{\phi}$ the system can be solved in the nodal DG-FEM framework. In $\Omega^{pml}$, Eq. \eqref{PML2norder} becomes:
\begin{subequations}
\begin{align}
    p_t &= -\rho c^2 \nabla \cdot \textbf v - \sigma_x p, \\
    \textbf v_t &= -\frac{1}{\rho} \nabla p - \frac{1}{\rho c^2} \boldsymbol{\phi}, \label{PMLvecV}\\
    \boldsymbol{\phi}_t &= \Gamma_1^x \boldsymbol{\phi} + c^2 \Gamma_2^x \nabla  p. \label{phiBefore}
\end{align}
\end{subequations}
Besides, to simplify the implementation of Eq. \eqref{phiBefore}, the calculation of $\nabla p$ is re-used from Eq. \eqref{PMLvecV} since it is equal to $-\rho \textbf v_t -\frac{1}{c^2} \boldsymbol{\phi}$. The auxiliary equation Eq. \eqref{phiBefore} becomes:
\begin{align}
    \boldsymbol{\phi}_t &= (\Gamma_1^x - \Gamma_2^x) \boldsymbol{\phi} - \rho c^2 \Gamma_2^x  \textbf v_t, \\
    \boldsymbol{\phi}_t &= -\begin{bmatrix}0&0&0\\0&\sigma_x&0\\0&0&\sigma_x\end{bmatrix} \boldsymbol{\phi} + \rho c^2 \begin{bmatrix}\sigma_x&0&0\\0&-\sigma_x&0\\0&0&-\sigma_x\end{bmatrix}  \textbf v_t.
\end{align}
The same procedure is repeated for a unique damping in the $y$-direction ($\sigma_y \neq 0$ and $\sigma_x = \sigma_z = 0$):
\begin{subequations} 
\begin{align}
    p_t &= -\rho c^2 \nabla \cdot \textbf v - \sigma_y p, \\
    \textbf v_t &= -\frac{1}{\rho}\nabla  p - \frac{1}{\rho c^2} \boldsymbol{\phi},\\
    \boldsymbol{\phi}_t &= -\begin{bmatrix}\sigma_y&0&0\\0&0&0\\0&0&\sigma_y\end{bmatrix} \boldsymbol{\phi} + \rho c^2 \begin{bmatrix}-\sigma_y&0&0\\0&\sigma_y&0\\0&0&-\sigma_y\end{bmatrix}  \textbf v_t,
\end{align}
\end{subequations}
and in the $z$-direction ($\sigma_z \neq 0$ and $\sigma_x = \sigma_y = 0$):
\begin{subequations}
\begin{align}
    p_t &= -\rho c^2 \nabla \cdot \textbf v - \sigma_z p, \\
    \textbf v_t &= -\frac{1}{\rho}\nabla  p - \frac{1}{\rho c^2} \boldsymbol{\phi},\\
     \boldsymbol{\phi}_t &= -\begin{bmatrix}\sigma_z&0&0\\0&\sigma_z&0\\0&0&0\end{bmatrix} \boldsymbol{\phi} + \rho c^2 \begin{bmatrix}-\sigma_z&0&0\\0&-\sigma_z&0\\0&0&\sigma_z\end{bmatrix}  \textbf v_t.
\end{align}
\end{subequations}
By linearly combining the three independent PML formulations with each other, the final PML formulation is:
\begin{subequations}
\begin{align}
    p_t &= -\rho c^2 \nabla \cdot \textbf v - (\sigma_x+\sigma_y+\sigma_z) p, \\
    \textbf v_t &= -\frac{1}{\rho}\nabla  p - \frac{1}{\rho c^2} \boldsymbol{\phi},\\
    \boldsymbol{\phi}_t &= -\Gamma_1 \boldsymbol{\phi} + \rho c^2 \Gamma_2 \textbf v_t, \label{phiAfter}
\end{align}
where 
\begin{equation}
    \Gamma_1 = \begin{bmatrix} \sigma_y+\sigma_z&0&0\\ 0&\sigma_x+\sigma_z&0\\ 0&0&\sigma_x+\sigma_y\end{bmatrix} , \quad  \Gamma_2 = \begin{bmatrix} \sigma_x-\sigma_y-\sigma_z&0&0\\ 0&\sigma_y-\sigma_x-\sigma_z&0\\ 0&0&\sigma_z-\sigma_x-\sigma_y\end{bmatrix}.
\end{equation}

\end{subequations}

\subsubsection{Discrete PML}
The PML terms in the governing equations and the auxiliary PDEs Eq. \eqref{phiAfter} must be integrated into the spatial and temporal schemes of the governing equation laid out in sections \ref{sec:SpatialDiscrete} and \ref{sec:TemporalDiscrete}. 
In a discrete context, $\sigma_i$ take up the value of the function $\sigma_i(x_i)$ at each node $x_{(k,n)}$ for $k\in \{1,...,K\}$ and $n\in \{1,...,N_p\}$.
At the initial time, the auxiliary variables are set to $0$. At each start of a new time step, the initial stage is equal to the results of the previous time step $t_n$:
\begin{subequations}
\begin{align}
    p^{(0)}_h &= p^n_h, \\ \textbf v^{(0)}_h &= \textbf v^n_h, \\ \boldsymbol{\phi}^{(0)}_h &= \boldsymbol{\phi}^n_h.
\end{align}
\end{subequations}
At the five stage $\alpha\in\{1,...,5\}$, the right-hand side of the acoustic governing equations are compute first in order to re-use the calculation of $\textbf v^{(\alpha)}_h$ in Eq. \eqref{Vtime} in the calculation of $k_\phi^{(\alpha)}$ for determining $\boldsymbol{\phi}^{(\alpha)}_h$:
\begin{subequations}
\begin{align}
    k_p^{(\alpha)} &= a_\alpha k_p^{(\alpha-1)} + \Delta t \left( \mathcal{L}_h(p_h^{(\alpha-1)},t_n) - (\sigma_x+\sigma_y+\sigma_z) p_h^{(\alpha-1)} \right), \\
    k_v^{(\alpha)} &= a_\alpha k_v^{(\alpha-1)} + \Delta t \left( \mathcal{L}_h(\textbf v_h^{(\alpha-1)},t_n) - \frac{1}{\rho c^2}\boldsymbol{\phi}_h^{(\alpha-1)} \right), \\
    p^{(\alpha)}_h &= p_h^{(\alpha-1)} + b_\alpha k_p^{(\alpha)}, \label{Ptime}\\
    \textbf v^{(\alpha)}_h &= \textbf v_h^{(\alpha-1)} + b_\alpha k_v^{(\alpha)}, \label{Vtime}\\
    k_\phi^{(\alpha)} &= a_\alpha k_\phi^{(\alpha-1)} + \Delta t \left( -\Gamma_1 \boldsymbol{\phi}^{(\alpha-1)}_h + \rho c^2 \Gamma_2 \textbf v^{(\alpha)}_h \right), \label{phiTime} \\
     \boldsymbol{\phi}^{(\alpha)}_h &= \boldsymbol{\phi}_h^{(\alpha-1)} + b_\alpha k_\phi^{(\alpha)}.
\end{align}
\end{subequations}

Finally, the fifth stage is the numerical solution at the new time step $t_{n+1}$:
\begin{subequations}
\begin{align}
    p_h^{n+1} &=p_h^{(5)}, \\ \textbf v_h^{n+1} &=\textbf v_h^{(5)}, \\ \boldsymbol{\phi}_h^{n+1} &=\boldsymbol{\phi}_h^{(5)}.
\end{align}
\end{subequations}

\subsubsection{Damping in the discrete PML} \label{sec:dampingPML}
Once the problem is discretized, the continuous damping functions $\sigma_x(x)$, $\sigma_y(y)$ and $\sigma_z(z)$ take discrete values at each nodes. In this section, constraints on the shape of the profile are determined by looking into the continuity of the solution at the interface of the PML zone and computational domain \cite{Engsig2006}. 

Let's consider the dampened solution variable $\textbf u^*(\textbf x_i)$ for  $\textbf x_i \in \Omega_{pml}$, it can be described as:
\begin{equation}
    \textbf u^*(\textbf x_i) = \zeta(\textbf x_i) \, \textbf  u(\textbf x_i),
\end{equation}
with $\zeta(\textbf x)$ a function in $[0,1]$. If the problem is considered in 1D, the dampened solution can be approximated in the modal form at time $t$ as:
\begin{equation} \label{relaxedModal}
    \textbf u^*(x_i) = \sum^{N_p}_{n=0}  \hat {\textbf u}_n(t) \psi_n(x_i) \zeta(x_i),
\end{equation}
where $x_i$ is the modal position in the 1D space, $\psi_n(x_i)$ is the local polynomial basis and $\hat u_n(t)$ is the expansion coefficient. To ensure the continuity of the solution variable $\textbf u$ at the interface between the domain of interest and the PML domain, the dampened solution $\textbf u^*$ and its spatial derivatives in both domains should be equals at all time $t$ so:
\begin{align} 
    \textbf u^*(x_i,t) &= \textbf u(x_i,t), \label{continuityU1}\\ 
    \frac{\partial \textbf u^*(x_i,t)}{\partial x} &= \frac{\partial \textbf u(x_i,t)}{\partial x}, \label{continuityU2}
\end{align}
As a consequence of the first continuity condition in Eq. \eqref{continuityU1} and the expression of the solution as a dampened solution Eq. \eqref{relaxedModal}, the following condition at the interface between the computational domain of interest and the PML zone applies:
\begin{equation} \label{zeta0}
    \zeta(x_i) = 1.
\end{equation}
Using the chain rule, the second continuity condition Eq. \eqref{continuityU2} on the spatial derivative of the dampened solution can be determined:
\begin{equation}
    \sum^{N_p}_{n=0} \hat u_n(t) \left[\psi_n(x_i) \frac{\partial \zeta(x_i)}{\partial x}  + \frac{\partial \psi_n(x_i)}{\partial x}  \zeta(x_i)\right] =  \sum^{N_p}_{n=0} \hat u_n(t) \frac{\partial \psi_n(x_i)}{\partial x}.
\end{equation}
So the criteria to satisfy at the interface is:
\begin{equation}
    \psi_n(x_i) \frac{\partial \zeta(x_i)}{\partial x} + \frac{\partial \psi_n(x_i)}{\partial x} \zeta(x_i) = \frac{\partial \psi_n(x_i)}{\partial x}.
\end{equation}
At the interface between the domain of interest and PML, $\zeta(x_i)=1$ so the derivative is:
\begin{equation} \label{zeta1}
     \frac{\partial \zeta(x_i)}{\partial x} = 0.
\end{equation}
As a consequence, at the discrete level, the damping profile $\sigma_i$ must start from $0$ at the entrance of the PML to respect Eq. \eqref{zeta0} and the $x_i$-derivative of $\sigma_i(x_i)$ applied to the coordinates of the entrance of the PML must be equal to $0$ to respect Eq. \eqref{zeta1}.

Extending those results to the 3D domain, the damping functions $\sigma_i(\textbf x_i)$ are null inside the domain of interest (i.e. when $\textbf x_i \in \Omega$). When $\textbf x_i \in \Omega^{pml}$, $\sigma_i(\textbf x_i)$ gradually grows from $0$ at the entrance of the PML to its maximum damping coefficient $\sigma_{max}>0$ at the termination of the PML.
Besides, the damping function must increase smoothly to avoid introducing local error that would generate numerical errors in the PML that could spread to the domain of interest.

\section{Method}\label{sec:method}
A common configuration is used to assess the damping performance, accuracy and stability of the PML formulation. The following setup and parameters are used for all tests in the results section unless specified otherwise.
\subsection{Mesh} \label{sec:Mesh}
The mesh is generated with Gmsh \cite{Geuzaine2009Gmsh}. The domain of interest $\Omega$ is a cube centered in $(0,0,0)m$ with an edge length of $L = 5 \;m$. The dimension of the tetrahedral mesh and the characteristic length of the elements are chosen to fulfill the following conditions: (1) the number of wavelength per element is roughly of 9 at peak frequency for accurate results, (2) no energy in the elements of the PML at the initial time (3) minimal number of elements in the mesh for computational efficiency.
\begin{center}
\begin{figure} [!ht]
    \centering
    \includegraphics[width=0.5\textwidth]{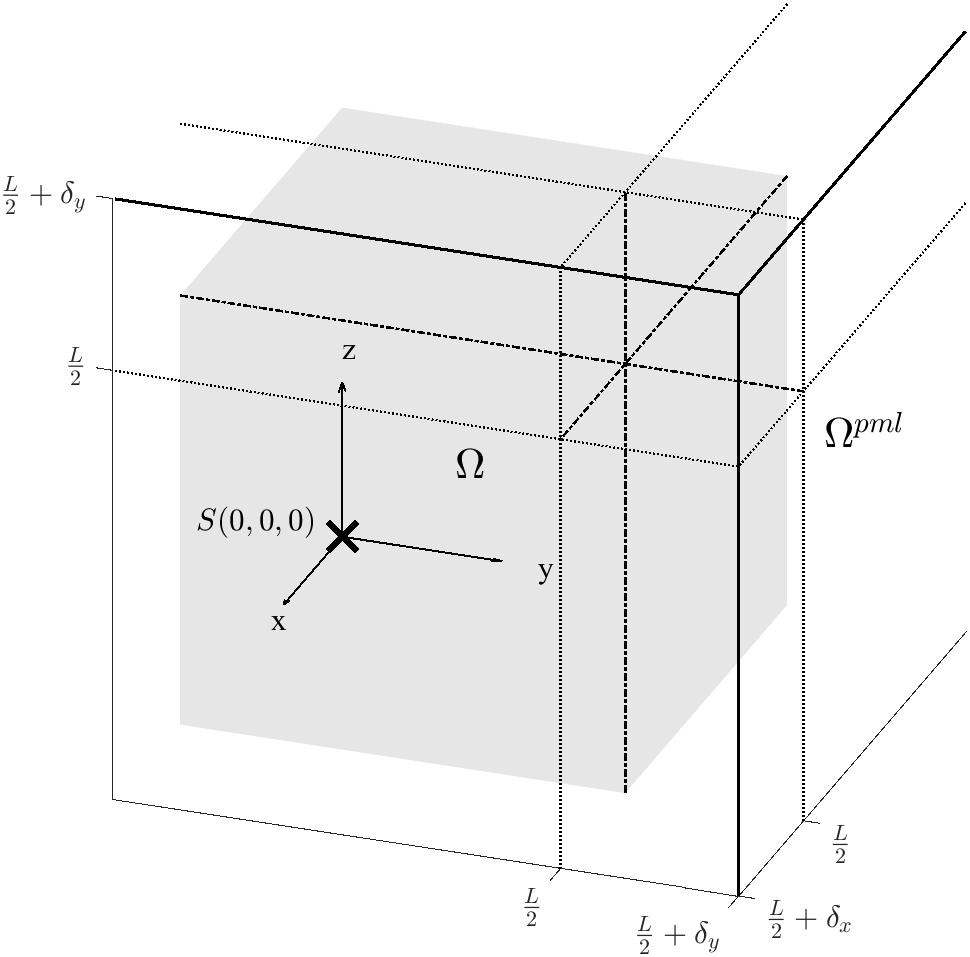}
    \caption{Notation in the mesh used. The computational domain (shaded) is surrounded by PML zones.}
    \label{fig:meshNotation}
\end{figure}
\end{center}

The PML's widths $\delta_x, \delta_y$ and $\delta_z$ are either equal to the peak wavelength of the source $\lambda_{peak}$ or half of $\lambda_{peak}$ and are the equal to one another in all three directions.

\subsection{Initial condition}  \label{sec:IC}
At initial time $t=0 \, s$, the source is a spatially distributed Gaussian pulse with a peak frequency $f_{peak}$ [$Hz$] where $\omega$ is the related pulse width. The pulse is centered at the source position $S = (0,0,0)m$ and the peak pressure at initial time is $P_0 = 1 \, Pa$ (equivalent to a sound pressure level of $94 \, dB$):
\begin{equation}\label{GaussianPulse}
p(\textbf x,0) = P_0 e^{ - \left(\frac{\pi \sqrt{2} f_{peak}}{c} \right)^2 r^2} = P_0 e^{ - \frac{1}{\omega^2} r^2},
\end{equation}

where $r = \sqrt{(x-x_s)^2+(y-y_s)^2+(z-z_s)^2}$ is the Euclidean distance between the source and the point at interest in the mesh. The particle velocity components in all three directions are null at initial time $t = 0$ s: $\textbf v (\textbf x,0) = \textbf 0 \, m/s$ and the auxiliary variables $\boldsymbol{\phi}$ introduced for the PML are also null at initial time: $\boldsymbol{\phi} (\textbf x,0) = \textbf 0$.

To vary the frequency content emitted from the source, three peak frequencies $f_{peak}$ are used: $343$, $514.5$ and $ 646 \, Hz$ with the corresponding peak wavelength $\lambda_{peak}$ equal to $1$, $\frac{1}{3}$ and $\frac{1}{2} \, m$ respectively.

\begin{center}
\begin{figure} [!ht]
    \centering
    \includegraphics[width=\textwidth]{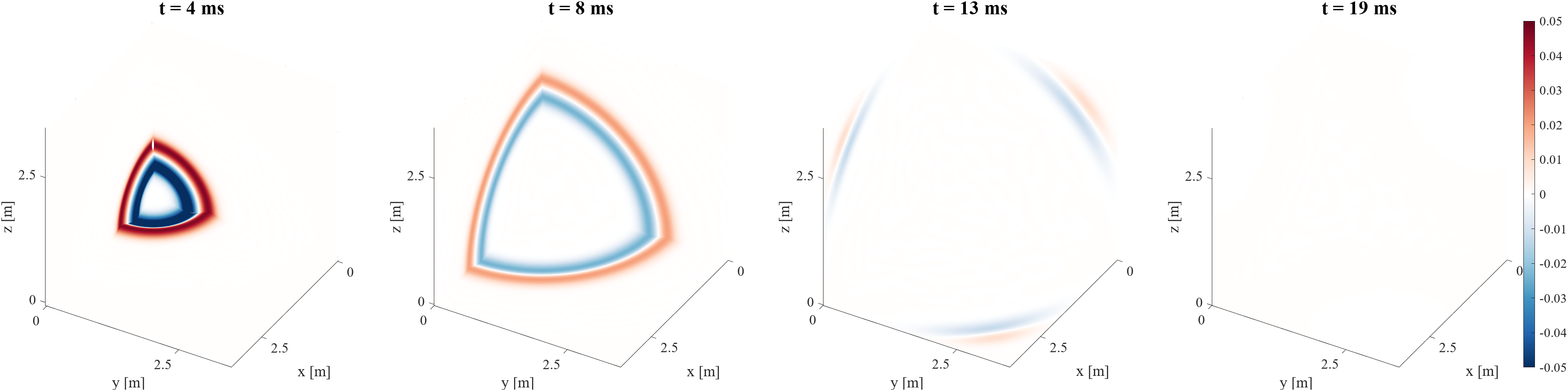}
    \caption{Snapshots of the time evolving 3D pressure field $p(\textbf x,t)$ in the positive 2D planes $x=0m$, $y=0m$ and $z=0m$ at different instants for an impulse generated at the source position $S$ and absorbed by a PML of width $1 \, m$ in all three directions.}
    \label{fig:pulseInTime}
\end{figure}
\end{center}

\subsection{Damping function $\sigma_i$} \label{sec:DampingFunction}
\begin{center}
\begin{figure} [!ht]
    \centering
    \includegraphics[width=0.5\textwidth]{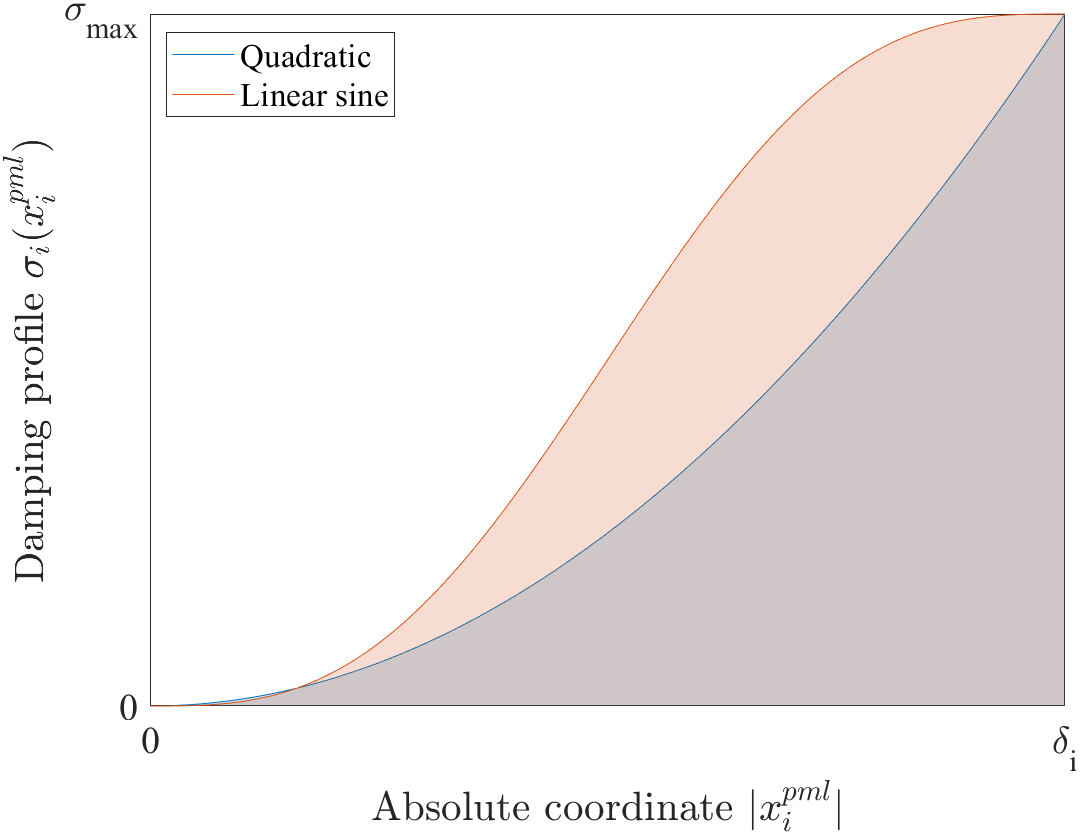}
    \caption{Shape of the two damping profiles tested in the local coordinates of the PML region.}
    \label{fig:dampingProfile}
\end{figure}
\end{center}

A judicious choice for the damping profile and its maximum damping coefficient can improve the effectiveness of the PML at no additional computational cost. Considering the continuity and smoothness conditions defined in Section \ref{sec:dampingPML}, two damping profiles are tested:
\begin{itemize}
    \setlength\itemsep{0em}
    \item quadratic: $\sigma_i(x^{pml}_i) = \sigma_{max} \left( \frac{|x^{pml}_i|}{\delta_i} \right)^2$,
    \item linear sine proposed in \cite{Grote2010}: $\sigma_i(x^{pml}_i) = \sigma_{max} \left( \frac{|x^{pml}_i|}{\delta_i} - \frac{sin \left(2\pi\frac{|x^{pml}_i|}{\delta_i} \right)}{2\pi}\right)$.
\end{itemize}
where $x^{pml}_i$ is the local coordinate inside the PML in each corresponding direction ($x_i= \{x,y,z\}$). $x^{pml}_i$ and  $\sigma_i(x^{pml}_i)$ start at zero  at the entrance of the PML and are equal to $\delta_{x_i}$ and $\sigma_{max}$ respectively at the termination of the layer, as illustrated in Figure \ref{fig:dampingProfile}.
For the choice of the maximum damping coefficient $\sigma_{max}$, the expression of the reflected wave arriving at normal incidence towards the PML as in \cite{Kaltenbacher2013} is considered. This wave is the one that has been the least dissipated since it has propagate along the shortest path in the domain to come back to its initial position. The reflected wave expression is as follows: 
\begin{equation}
    P_R = P_0e^{-2/c \int^{\delta_i}_0\sigma_i(x^{pml}_i)dx^{pml}_i} = P_0R,
\end{equation}
where R is defined as the reflection factor.

By taking the damping as constant ($\sigma_i(x^{pml}_i) = \sigma_0$), a reference maximum damping coefficient $\sigma_0$ can be extracted: $\sigma_0=-\frac{c}{2}\frac{ln(R)}{\delta_i}$ with $R$ being $10^{-3}$. In the results section, the maximum damping coefficient $\sigma_{max}$ is calibrated with $\sigma_0$.

To normalized the impact of the damping on the wave the damping "area" is defined. It is the region underneath the damping profiles $\sigma_i(x^{pml}_i)$ so it can be calculated as:
\begin{equation}
    \int^{\delta_i}_0\sigma_i(x^{pml}_i) dx^{pml}_i = \begin{cases} \frac{1}{3} \sigma_{max}\delta_i \text{ for a quadratic profile,} \\ \frac{1}{2} \sigma_{max}\delta_i \text{ for a linear sine profile.}
    \end{cases}
\end{equation}

\section{Results} \label{sec:Results}
The PML formulation is assessed in terms of (1) damping configurations to optimize the damping function, (2) convergence and accuracy, (3) ability to absorb grazing incidence wave, and (4) long-time stability after a large amount of time steps for different applications.
\subsection{Optimization of the damping function} \label{sec:dampOpti}
When a PML dampens outgoing waves, the energy carried by the waves vanishes in the PML domain. Therefore, the instantaneous energy calculated in the domain of interest should decrease rapidly as the waves reach the PML. The energy remaining has two origins: the energy related to the numerical error inherent to the discretization and the energy related to the PML absorption. The present subsection seeks to minimise the second one.

The total instantaneous acoustic energy $E(t)$ in the domain of interest at time $t$ is the sum of the instantaneous potential and kinetic energy at the $N_p$ nodes in the $K$ elements belonging to the domain of interest $\Omega$ and is defined as:
\begin{equation}\label{EnergyDG}
    E(t) = \sum^{K}_{k=1} \sum^{N_p}_{n=1} \frac{1}{2} \frac{1}{\rho c^2} p(\textbf x_{(k,n)},t)^2 +  \frac{1}{2} \rho \lVert \textbf v(\textbf x_{(k,n)},t) \rVert ^2 .
\end{equation}

\begin{center}
\begin{figure} [htbp]
    \centering
    \includegraphics[width=0.95\textwidth]{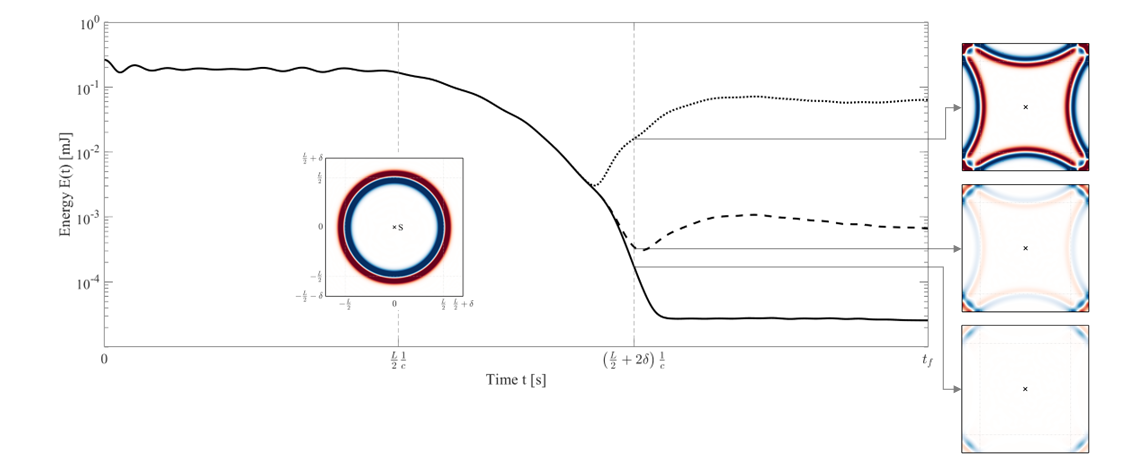}
    \caption{Energy $E(t)$ in $\Omega$  for three cases: (dotted line) no PML treatment, (dashed line) an intermediary cases and (solid line) optimal PML damping function. For all cases, the termination of the PML region is perfectly reflective and the source is in position $S(0,0,0)$. Snapshots of the pressure $p(\textbf{x},t)$ in a 2D plane illustrate the state of the wave when it first reaches the PML ($t=\frac{L}{2} \frac{1}{c}$) and the three cases when the normal wave has bounced back and crosses the PML interface with the inner domain a second time ($t=\left( \frac{L}{2}+2\delta \right) \frac{1}{c}$).}
    \label{fig:xiR_explanation}
\end{figure}
\end{center}

\begin{center}
\begin{figure} [!ht]
    \centering
    \includegraphics[width=0.95\textwidth]{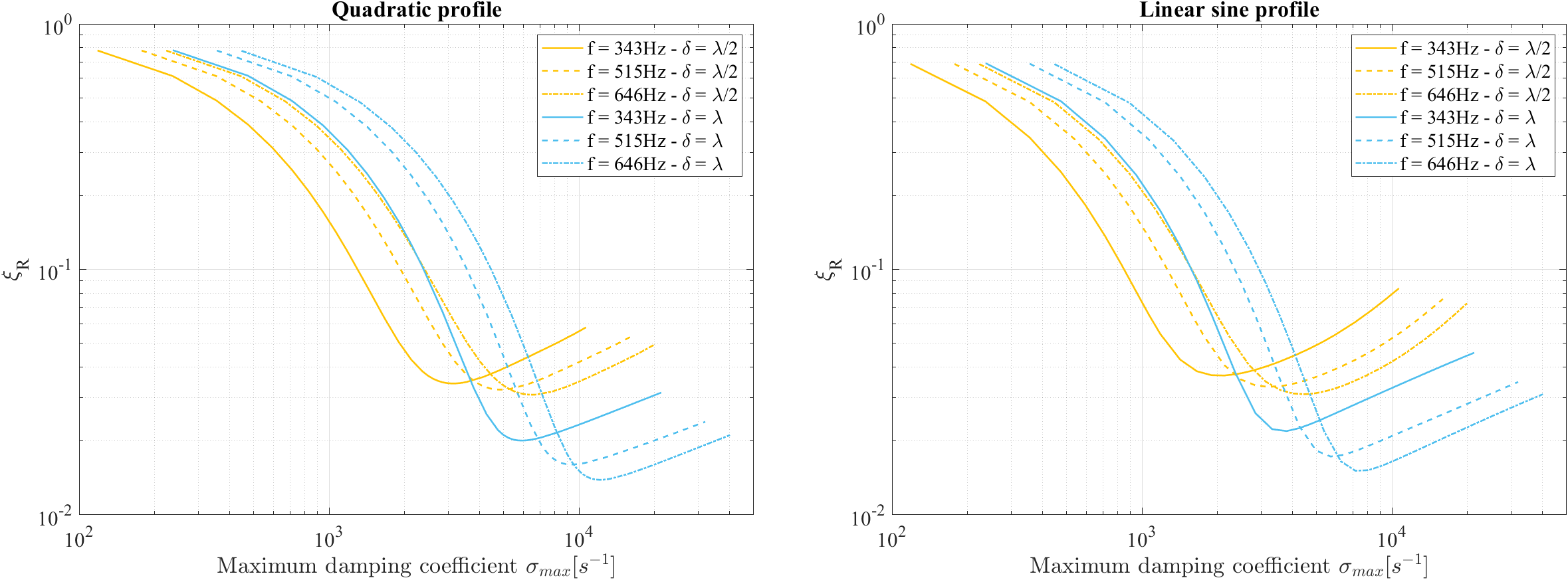}
    \caption{Relative error $\xi_R$ as a function of the maximum damping coefficient for two damping profiles and two width of PML (one wavelength or half a wavelength of the peak frequency) for three different peak frequencies.}
    \label{fig:dampingCoefficient}
\end{figure}
\end{center}

\begin{center}
\begin{figure}[!ht]
    \centering
    \includegraphics[width=0.95\textwidth]{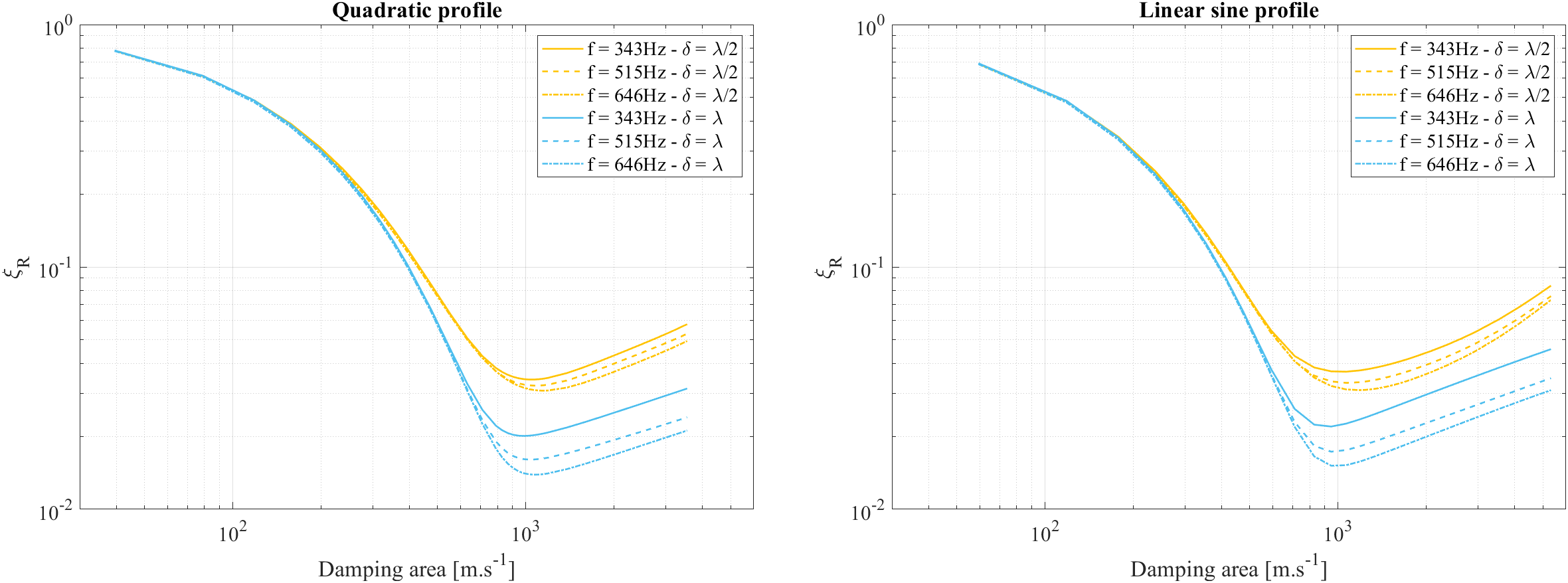}
    \caption{Relative error $\xi_R$ from Figure \ref{fig:dampingCoefficient} as a function of the damping area for two damping profiles.}
    \label{fig:dampingArea}
\end{figure}
\end{center}

The energy related to the PML absorption that remains in $\Omega$ after the dissipation of the wave by the PML varies depending on the damping profile and the damping maximum coefficient. As in \cite{Modave2014}, the boundary condition at the termination of the PML is perfectly reflective in order to solely account for the absorption caused by the PML. To quantify the effectiveness of the PML, we also estimated the relative error $\xi_R$ defined as:
\begin{equation}
    \xi_R=\sqrt{\frac{E(t_f)}{E_{wall}(t_f)}},
\end{equation}
where $t_f$ is the time at which the wavefront arriving normally towards the PML is reflected back to the initial position so at time $t_f=(L+2\delta_i)/c$ as shown in Figure \ref{fig:xiR_explanation}. $E_{wall}$ is the total energy in $\Omega$  with no PML treatment of the layer and, unlike \cite{Modave2014}, a reflective boundary condition is imposed at the termination of the PML (i.e. dotted line case in Figure \ref{fig:xiR_explanation}).
$\xi_R$ measures the capacity of the PML to absorb the wave so $\xi_R$ should tend to decrease with increasing $\sigma_{max}$.
The relative error $\xi_R$ is computed for $N=3$ and for a range of maximum damping coefficient from $0.1$ to $10$ times the reference maximum damping coefficient $\sigma_0$, determined in Section \ref{sec:DampingFunction}. The PML width is equal to one wavelength (in blue in Figure \ref{fig:dampingCoefficient} and Figure \ref{fig:dampingArea}), corresponding to 4 layers of elements in the PML or equal to half of the wavelength  (in yellow in Figure \ref{fig:dampingCoefficient} and Figure \ref{fig:dampingArea}), corresponding to 2 layers of elements in the PML.
In Figure \ref{fig:dampingCoefficient}, the relative error $\xi_R$ logically tends to 1 when the maximum damping coefficient $\sigma_{max}$ tends to 0 since the PML has little impact on the wave when the damping is negligible. As was expected, the error decreases rapidly as the damping maximum increases until reaching an optimal value. Then, $\xi_R$ increases at different rate depending on the damping profile. The degradation of the error after a certain maximum damping coefficient can be explained by the limitation of the discretization: when the wave decays in the PML, the jump of damping inside the layer can be too strong and non-negligible numerical errors are introduced. For each damping profile and width of PML, an optimal value of $\sigma_{max}$ can be defined. 

As noted in Section \ref{sec:DampingFunction}, the area underneath the damping profile differs depending on the damping profile and the maximum damping coefficient $\sigma_{max}$, itself dependent on the peak frequency since the PML widths are chosen as a multiple of the peak wavelength. Figure \ref{fig:dampingArea} re-scales the curves of Figure \ref{fig:dampingCoefficient} to draw the relative error $\xi_R$ as a function to the respective damping area.

The cases with a smaller PML width (in yellow in Figure \ref{fig:dampingCoefficient} and \ref{fig:dampingArea}) has an error almost twice as big as the equivalent case with the a wider PML (in blue in the same figures). The width of the PML is therefore determining factor on the efficiency of the PML and should be chosen according to the accuracy tolerated for the simulation and the computational capabilities.

Looking at the error $\xi_R$ from a damping area point-of-view, a common damping optimum can be observed at around $1 \, 000 \; m/s$ in Figure \ref{fig:dampingArea}. Both profiles (quadratic and linear sine) present comparable error so this suggests that the value of the maximum damping coefficient has a greater impact on the performance of the PML than the damping profile. 

The quadratic profile is slightly more efficient than the linear sine profile so in the following sections, the damping function used in the PML have a quadratic profile and the maximum damping coefficient is the one minimizing the error $\xi_R$ from an damping area point-of-view.

\subsection{Convergence} \label{sec:Convergence}
In the continuous domain, PMLs are designed to absorb completely the \textit{exact} solution of the wave equation. However, when the problem is discretized with DG-FEM or other numerical methods, a discretization error is inevitably introduced and the wave equation in \textit{approximately} solved. As a result, the PML no longer perfectly absorb the outgoing waves. In this section, the convergence and computation time of the PML are studied in terms of the order of approximation.

For a fixed number of elements $K$ (i.e. for a fixed mesh) in a cube mesh described in subsection \ref{sec:Mesh}, the polynomial order of approximation $N$ is increased ($p$-refinement) to enhance the quality of the approximation of the solution. For all simulations, the peak frequency of the Gaussian pulse of the source is $343 \, Hz$, the PML width is equal to the corresponding peak wavelength: $1 \, m$ and absorbing boundary condition are imposed at the termination of the PML. Table \ref{tab:meshesConv} summarizes the properties of the four meshes tested in this convergence study. 
\begin{table}[!ht]
    \centering
    \begin{tabular}{|c|c|c|c|c|} \hline 
                                  & $M_1$& $M_2$& $M_3$& $M_4$\\ \hline 
 Minimal characteristic length [m]& 0.030& 0.023& 0.017& 0.016\\ \hline 
 Average characteristic length [m]& 0.095& 0.073& 0.064& 0.053\\ \hline 
 Maximum characteristic length [m]& 0.175& 0.123& 0.106& 0.090\\ \hline 
 Minimum wavelengths per element for $f_{peak}=343Hz$&  5.70&  8.13&  9.41&  11.1\\ \hline  Number of elements $K$ in $\Omega$& 4638& 10244& 15779& 26943\\ \hline\end{tabular}
    \caption{Characteristic length of the elements and number of elements in the domain of interest for the four meshes tested}
    \label{tab:meshesConv}
\end{table}

The computed solutions in the four meshes of Table \ref{tab:meshesConv} are compared to the corresponding computed solution in a reference mesh (subscript ''$ref$'') containing the exact same meshing of the domain of interest $\Omega$. The reference mesh has a $4.5\,m$-thick surrounding region with no PML treatment and only absorbing boundary condition at the termination of it, instead of a $1\,m$-thick PML. At time $t_f$ defined in the subsection \ref{sec:dampOpti}, the wave reaches the outer boundary of the reference mesh so the reference pressure $p_{ref}$ is not polluted by reflections from non-reflective treatments. 

\begin{center}
\begin{figure} [!htbp]
    \centering
    \includegraphics[width=0.75\textwidth]{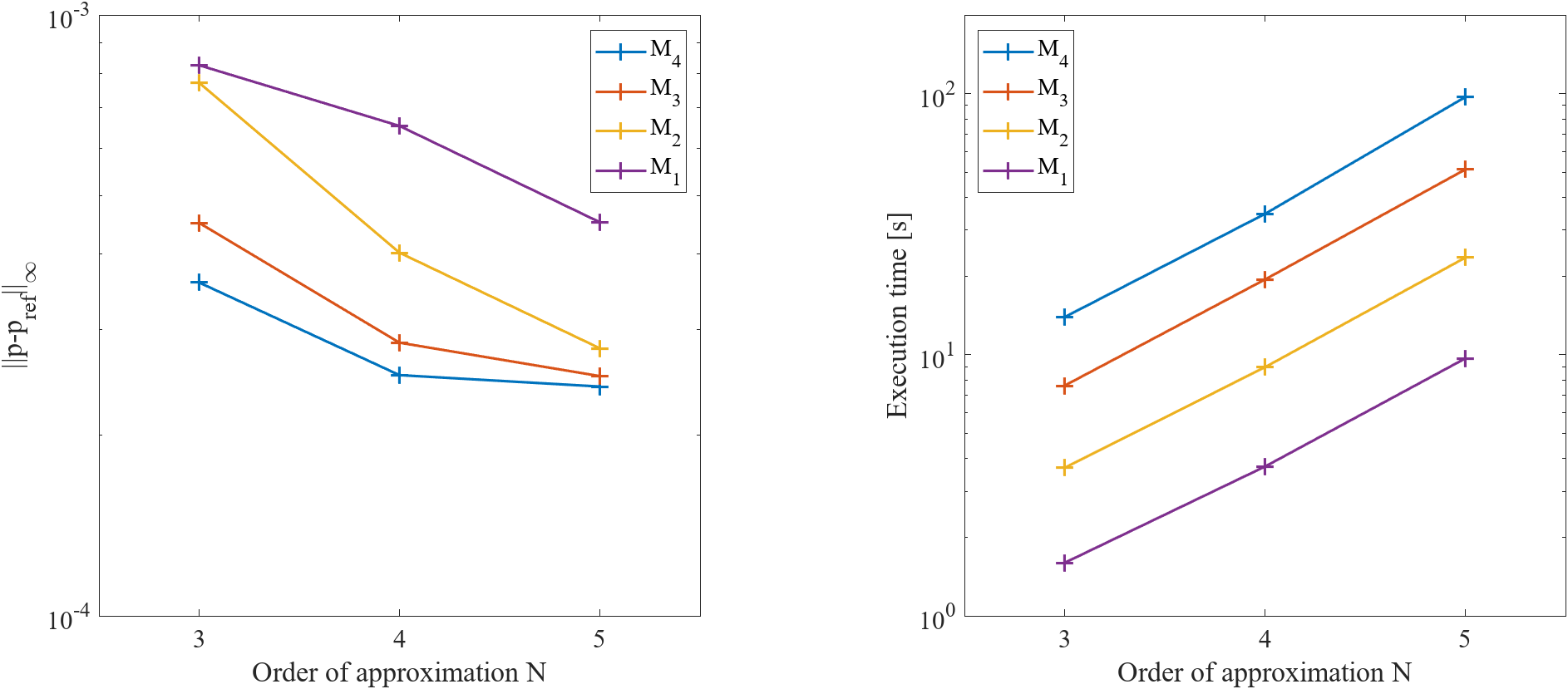}
    \caption{$L^\infty$-norm of the acoustic pressure in the domain of interest for different order of approximation and different meshes and the corresponding time of execution of the simulation for a final time of $t_f$.}
    \label{fig:pRefinement}
\end{figure}
\end{center}

The convergence error is defined as the $L^\infty$-norm of the difference between the acoustic pressure in the domain of interest and the acoustic pressure in the reference mesh at the final time of the simulation $t_f$: 
\begin{equation}
    L^{\infty}_P=\lVert p-p_{ref}\rVert_\infty = \max_{\textbf{x}\in \Omega} \left(|p(\textbf{x},t_f)-p_{ref}(\textbf{x},t_f)| \right).
\end{equation}
Figure \ref{fig:pRefinement} presents the $p$-refinement of the four meshes and the corresponding execution time. The simulations are run on a NVIDIA GeForce RTX 3070 Laptop GPU with 8GB dedicated memory. The solver is based on libParanumal \cite{libParanumal} and take advantage of the fact that DG-FEM is a highly parallelizable formulation. 

Figure \ref{fig:pRefinement} confirms that for the same mesh, a more accurate approximation of the wave (higher $N$) results in a reduced $L^\infty$-error on the acoustic pressure $p(\textbf{x},t)$ at time $t_f$. However increasing the order of approximation $N$ also increases exponentially the execution time it takes to complete the simulation. The rate of $p$-refinement is lower for finer mesh, meaning that when the mesh is fine, increasing the order of approximation will not improve the accuracy significantly enough to justify the increase in computational time.

\subsection{Grazing incidence}
When a wave approaches the outer boundary of the domain at a grazing incidence, the absorption in the direction of the propagation of the wave is inefficient and reflections from the termination of the PML can be generated, polluting the inner domain. To highlight the absorbing properties of the PML formulation for waves approaching the outer boundary of the domain at grazing incidence, the initial mesh is elongated along one of the axis. 

\begin{center}
\begin{figure} [htbp]
    \centering
    \includegraphics[width=\textwidth]{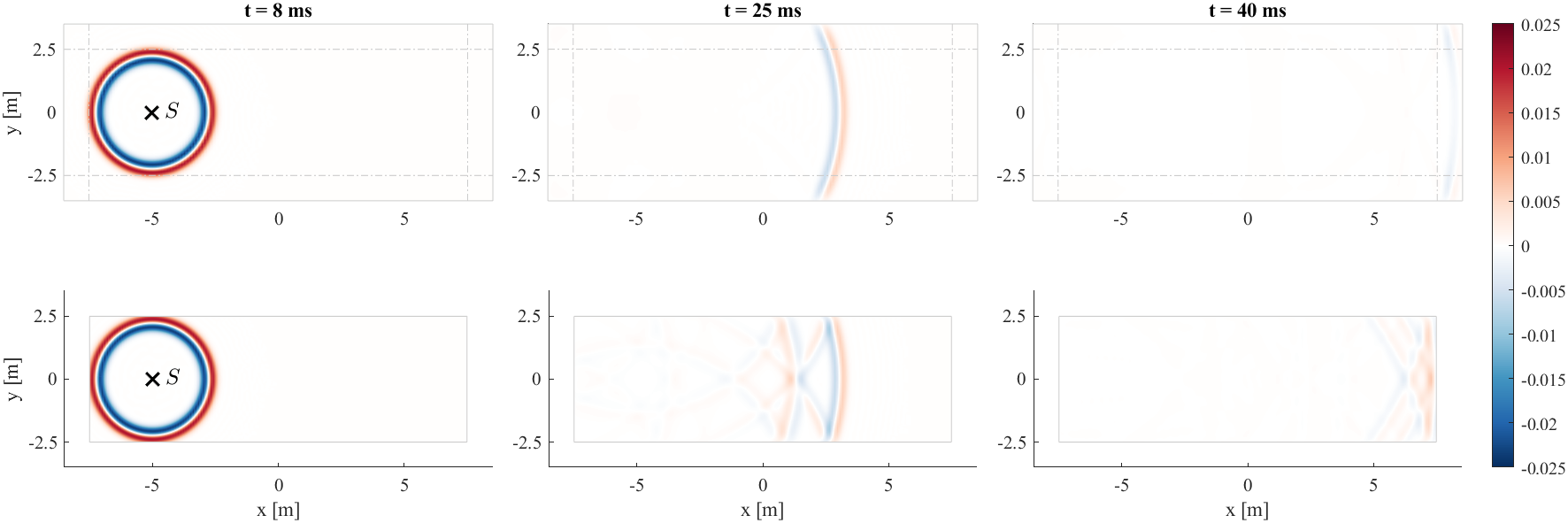}
    \caption{Snapshots of the acoustic pressure $p(\textbf x,t)$ at three stages of the simulations in the $z = 0m$ plane for the case with PML and ABC at its termination (top) and for the case with only ABC (bottom).}
    \label{fig:grazingPres}
\end{figure}
\end{center}

The elongated mesh is larger along the $x$-axis so the domain of interest becomes $15 \; m$-long instead of $L = 5 \; m$-long. The source position is shifted $5 \, m$ so that for the majority of its propagation, the wavefront is grazing the boundary of the domain along the x-axis. As in the Section \ref{sec:Convergence}, the peak frequency of the Gaussian pulse of the source is $343 \, Hz$, the PML width is equal to the corresponding peak wavelength: $1 \, m$ and absorbing boundary condition are imposed at the termination of the PML. Figure \ref{fig:grazingPres} shows the interpolated pressure on a 2D slice of the domain of the plane at $z=0$ m. 


\begin{center}
\begin{figure} [htbp]
    \centering
    \includegraphics[width=0.75\textwidth]{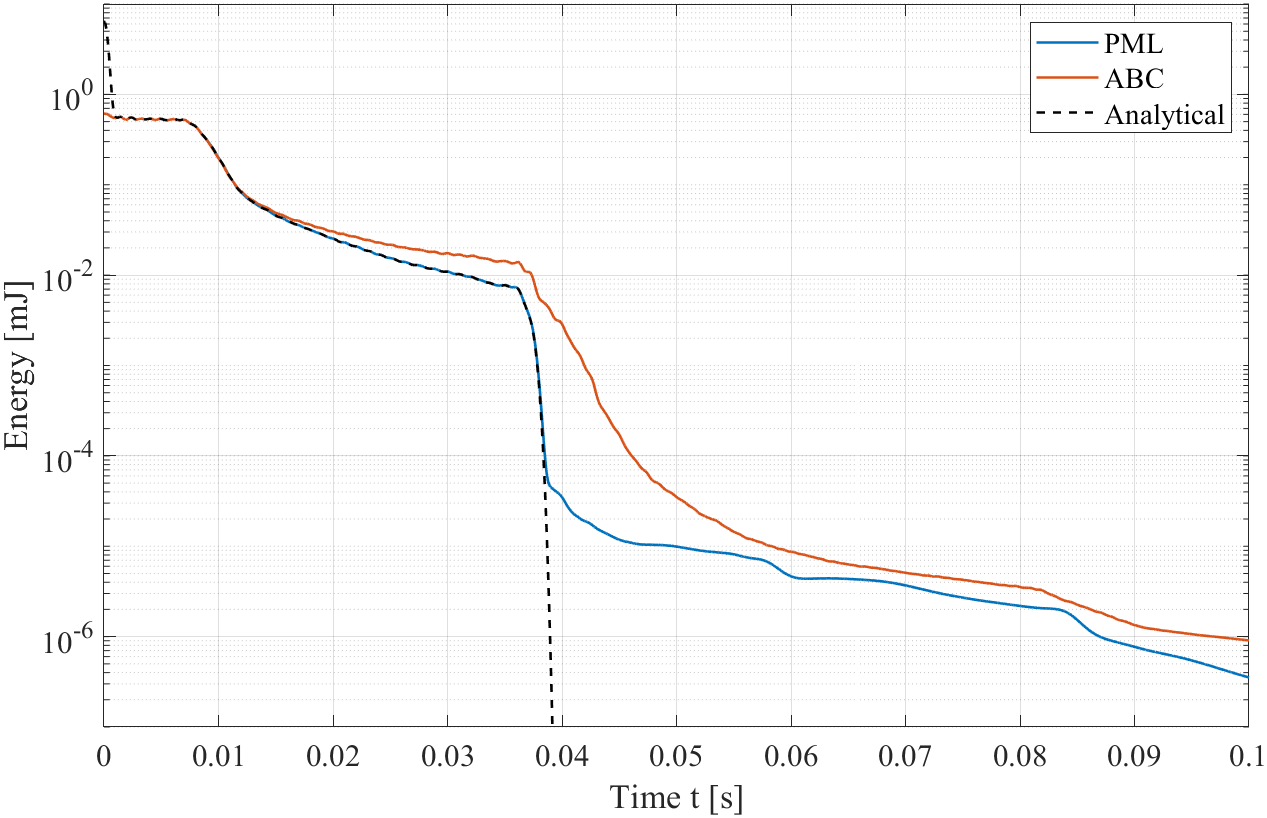}
    \caption{Energy inside the domain of interest over $0.1 \; s$ for an elongated mesh with a PML and with ABC. In dashed black line is the analytical energy calculated at the same nodes of the domain of interest.}
    \label{fig:grazingE}
\end{figure}
\end{center}

The energy $E(t)$ of the numerical solution in the elongated mesh is compared to energy of the exact solution, calculated at the exact same nodes position and time steps. For a source at position $S(x_s,y_x,z_s)$ and a Gaussian distribution of the acoustic pressure at initial time of width $\omega$, the analytical pressure and velocity are calculated as follows:
\begin{subequations}
\begin{align}
    p_{exact}(\textbf x,t) &= \frac{1}{2r}\left( (r-ct)e^{-\frac{(r-ct)^2}{\omega^2}}+(r+ct)e^{-\frac{(r+ct)^2}{\omega^2}}\right), \\
    \textbf v_{exact}(\textbf x,t) &= \frac{1}{2rc\rho} \left( \left(\frac{\omega^2}{2r} + r-ct\right)e^{-\frac{(r-ct)^2}{\omega^2}} + \left(\frac{\omega^2}{2r} + r+ct\right) e^{-\frac{(r+ct)^2}{\omega^2}}\right) \frac{\textbf x - \textbf x_s}{r},
\end{align}
\end{subequations}
where $r = \sqrt{(x-x_s)^2+(y-y_s)^2+(z-z_s)^2}$ is the Euclidean distance between the source and the point at interest in the mesh. The energy $E(t)$ over $0.1\,s$ and for $N=4$ is shown in Figure \ref{fig:grazingE}.

In Figure \ref{fig:grazingE}, $E(t)$ dips first at around $0.01 \; s$ when a portion of the wave has already reached the closest PML, i.e. the PML in the $y$ and $z$-directions and the PML in the negative $x$-direction. The second dip in energy at around $0.038 \; s$ translates the absorption of the wave at the end of the elongated domain along the x-axis. The difference in energy between the PML case and the ABC can be explained by the large spurious reflection generated by the ABC as shown in the bottom Snapshots of Figure \ref{fig:grazingPres}. After the second energy dip, an unacceptable amount of energy still remains in the domain of interest for the case of the ABC which pollutes the computed solution if an outdoor acoustic simulation is performed. For the case of the PML, the energy curve is closely following the analytical curve so  a small amount of energy is reflected back into the domain when the waves propagates in the PML as illustrated in the top snapshots of Figure \ref{fig:grazingPres}. 
Because of the discretization and the fact that the PML can't absorb perfectly the grazing incidence wave, some error still remains after $0.04 \, s$ for the PML case. It can be lowered by refining the mesh or increasing the order of approximation but at a computational cost. The choice of PML over ABC simulation is therefore justified despite the extra computation if accurate outdoor acoustic simulation is needed.

\subsection{Long-time simulations}
Three simulations were performed with a simulation final time of $2\;s$ to check the long-time stability of the PML formulation which corresponds to a wave traveling $686 \,m$ away from the source position. For all three cases, the simulation is performed using the mesh $M_3$ from Section \ref{sec:Convergence} an an order of approximation of $N=3$ to have a reasonable computation time. The PML width is still $1\;m$  and absorbing boundary conditions are again imposed at the termination of the PML.

The first case is a Gaussian pulse excitation at initial time as in Eq. \eqref{GaussianPulse} with a peak frequency of the Gaussian pulse of $f_{peak} = 343 \, Hz$. Figure \ref{fig:longImpulse} presents the energy $E(t)$ inside the domain of interest $\Omega$ over the time of the simulation for a impulse at initial time over $2 \,s$ corresponding to $87\,970$ time steps:
\begin{center}
\begin{figure} [htbp]
    \centering
    \includegraphics[width=0.85\textwidth]{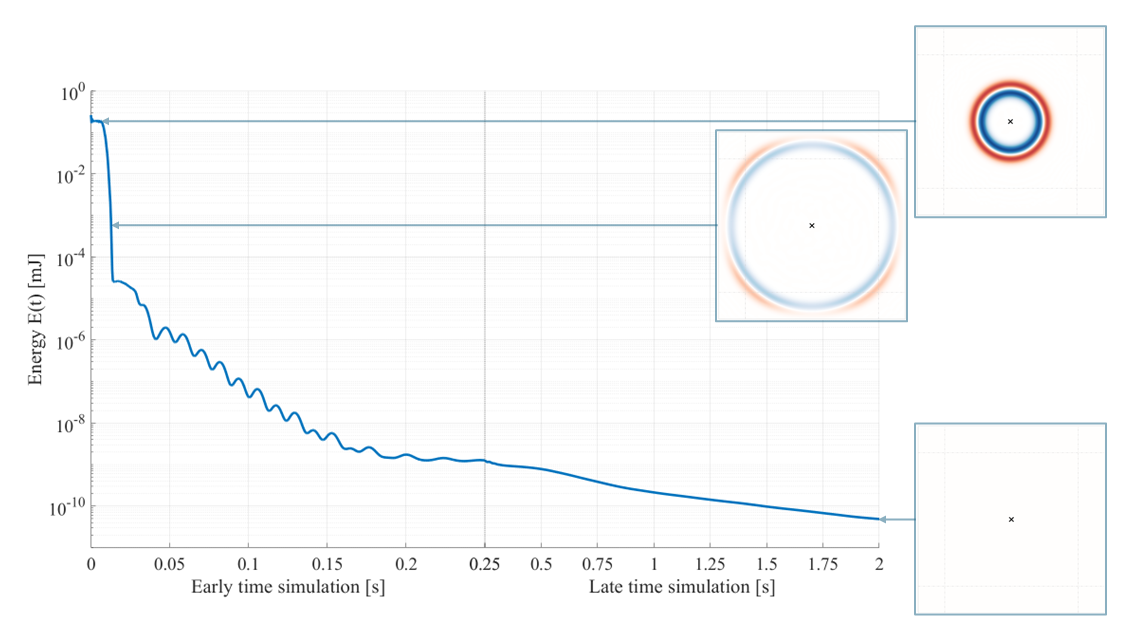}
\caption{Energy inside the domain of interest over $2\;s$ in the early stages of the simulation (left) and for the rest of the simulation time (right) for a impulse source. Snapshots of the acoustic pressure $p(\textbf{x},t)$ in a 2D plane illustrating different stages of the simulation}
 \label{fig:longImpulse}
\end{figure}
\end{center}

\begin{center}
\begin{figure} [htbp]
    \centering
    \includegraphics[width=0.85\textwidth]{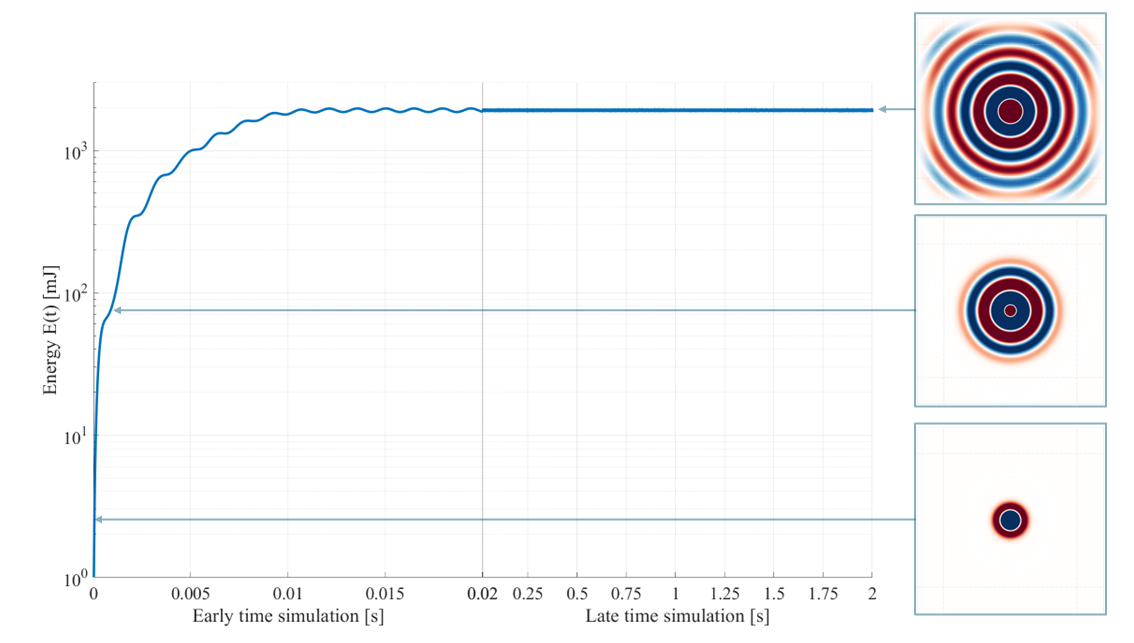}
    \caption{Energy inside the domain of interest over $2\;s$ in the early stages of the simulation (left) and for the rest of the simulation time (right) for a periodic source. Snapshots of the acoustic pressure $p(\textbf{x},t)$ in a 2D plane illustrating different stages of the simulation}
    \label{fig:longPulse}
\end{figure}
\end{center}

For the second case, the previous setup is kept the same but the source is changed a periodic Gaussian excitation by adding the following source term $S(\textbf x, t)$ to the acoustic pressure $p(\textbf{x},t)$ at each time step in the time loop:
\begin{equation}
    S(\textbf x, t) = p(\textbf x,0) \cos\left(\frac{1}{2\pi} f_{peak} t\right). 
\end{equation}
Figure \ref{fig:longPulse} presents the energy $E(t)$ inside the domain of interest $\Omega$ over the time of the simulation for a periodic pulse over $2 \,s$ also corresponding to $87\,970$ time steps.

\begin{center}
\begin{figure} [htbp]
    \centering
    \includegraphics[width=0.85\textwidth]{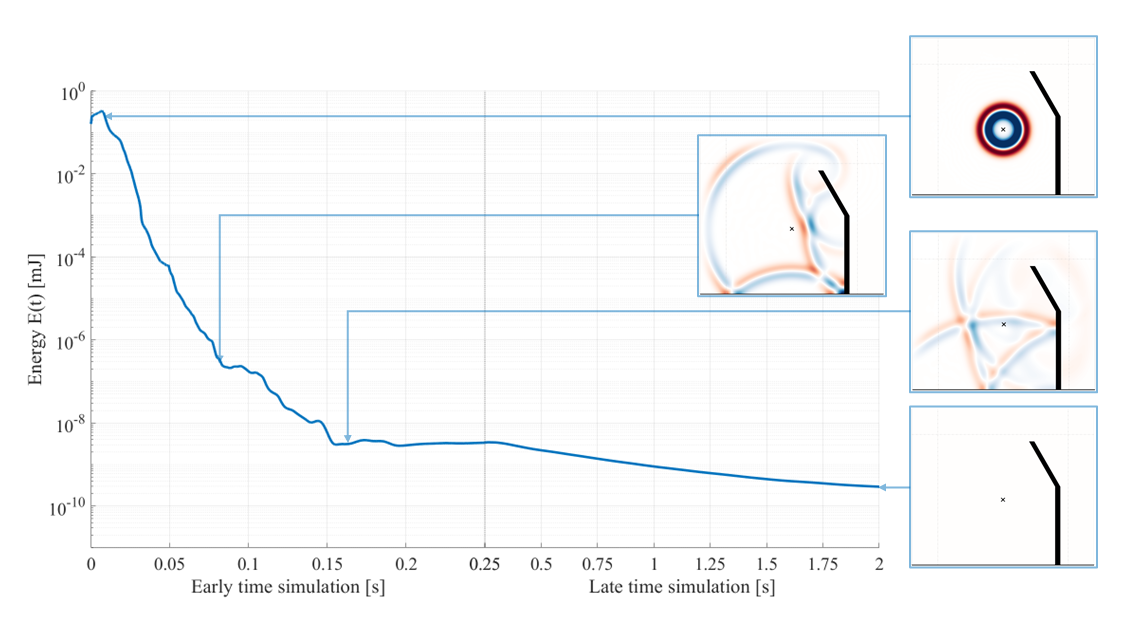}
    \caption{Energy inside the domain of interest over $2\;s$ in the early stages of the simulation (left) and for the rest of the simulation time (right) for an impulse source next to a sound barrier. Snapshots of the acoustic pressure $p(\textbf{x},t)$ in a 2D plane illustrating different stages of the simulation}
    \label{fig:longBarrier}
\end{figure}
\end{center}

One of the advantages of DG-FEM is the ability of this method to handle more complex geometry with high-order accuracy ; in particular geometries with non-right angles. This last example presents the simulation of an impulse source as in the first example next to a sound barrier $4 \, m$ long and $20 \, cm$ thick. It is tilted by $60^{\circ}$ angle towards the source to replicate a typical sound barrier shape. The ground and the barrier surface are perfectly reflecting and the mesh $M_3$ is modified accordingly.
Figure \ref{fig:longBarrier} presents the energy $E(t)$ inside the domain of interest $\Omega$ over $2\,s$ corresponding to $140\,317$ time steps. The increase in number of time step is explained by the necessity of having a finer mesh resolution around the sound barrier to conform the mesh to the shape of the barrier. Mesh refinement lowers the value of the minimum characteristic length of the tetrahedral elements which in turn decreases the time step value Eq. \eqref{timeStepDef} and a higher number of step are necessary to reach $2 \; s$ simulation time.

Figure \ref{fig:longImpulse} shows that the energy $E(t)$ inside the domain of interest trends downwards, suggesting that this PML formulation is stable for long-time simulations. For a periodic source in Figure \ref{fig:longPulse}, the decoupled PML effectively absorbs the outgoing radiating wave without introducing errors that could get out of control and pollute to the domain of interest. Indeed, the energy $E(t)$ is fluctuating periodically around $1\,900 \; mJ$ after $0.01 \, s$ and the PML formulation is stable.
In Figure \ref{fig:longBarrier}, the instantaneous energy takes a longer time to dissipate compared to in Figure \ref{fig:longImpulse} since the wave is reflected by the ground and the sound barrier (in thicker black line in Figure \ref{fig:longBarrier}). Nevertheless, $E(t)$ still trends downwards after most of the energy has left the domain of interest and the perfectly reflective boundaries do not affect the stability of the PML.

\section{Discussion}\label{sec:Discussion}
The decoupled PML formulation defined in Section \ref{sec:NumericalMethod} and examined in Section \ref{sec:Results} is stable and accurate but some limitations must be noted.

Section \ref{sec:dampOpti} defined an optimal function and maximum coefficient for the damping function by looking at the PML performance from a damping area point-of-view. The optimal value of the damping area is independent of the width of the PML and the frequency of the source but this optimization procedure is limited to the setup used. The value for the damping area is optimal only if the criteria chosen is to minimize the amount of energy coming back into the domain for a normal incidence wave. This study does not allow us to draw conclusion on the optimal damping function for cases where grazing incidence waves are predominant.

In terms of formulation, the decoupled PML differs from the PML of \cite{Grote2010} and \cite{Baffet2019} in the treatment of the edges and corners of the PML domain where the non-null damping terms are not coupled together. Indeed, since the PML formulation combines three stable PML defined independently in the three Cartesian directions, the cross terms of the original formulation \cite{Grote2010} are neglected which lowers the number of auxiliary variables, the number corresponding PDEs and number of PML terms to add to the governing equation in the PML zone.

The different PML formulation start from the same coordinate stretching Eq. \eqref{coordStretch} but the operations following the coordinate stretching generate different total number of auxiliary variables and corresponding PDEs. For the acoustic wave equation in 3D, a U-PML \cite{Hu2001,Hu2005,Parrish2009} or M-PML \cite{Meza2008} or a C-PML \cite{Komatitsch2007,Martin2009} create one auxiliary variable for each acoustics variable and for each of the three directions so 12 auxiliary variables and 12 corresponding PDEs are necessary. By comparison, our PML formulation has only 3 auxiliary variables ($\phi_x$, $\phi_y$ and $\phi_z$) and 3 corresponding PDEs Eq. \eqref{phiAfter}. Moreover, because of the simplicity of the formulation, the auxiliary PDEs can be calculated using anterior calculations inside the time loop, making the simulation even more computationally efficient.
To compare the efficiency in terms of execution time, if we simulate $8\,000$ time steps for a mesh of $183\,568$ tetrahedral elements, the simulation would be comparable to the application example of \cite{Modave2017} which lasted for $1\,h\,40\,min$ or $6\,000 \; s$. With the decoupled PML formulation and using the same hardware setup of Section \ref{sec:Convergence} and a highly parallelized solver, the execution time is $246 \; s$ for $N=3$, $469 \; s$ for $N=4$, and $957 \; s$ for $N=5$, which is significantly faster than \cite{Modave2017}.

\section{Conclusion}\label{sec:Conclusion}
We have proposed a decoupled PML formulation using DG-FEM for the acoustic wave equation in 3D that is stable for long-time simulation and a procedure to determine an optimal damping function to minimise the error generated by the PML for a wave approaching the layer normally. The PML formulation was evaluated in terms of damping function optimization, convergence, grazing incident wave absorption and long-time stability. By decoupling the damping terms in the corners and edges of the PML, the computational time is reduced without sacrificing accuracy. Taking advantage of the highly parallelizable framework and high-order accurate properties of the nodal DG-FEM, the outdoor acoustics solver is efficient and accurate which opens up opportunities for solving complex and large environmental acoustic problems that are usually limited in space and meshing conformity.

\section*{Acknowledgment}
The authors would like to thanks Treble Technologies for providing access to the code of their room acoustics DG-FEM solver.

\section*{CRediT authorship contribution statement}
\noindent \textbf{Sophia Julia Feriani}: Conceptualization, Methodology, Software, Validation, Investigation, Writing - Original Draft, Writing - Review \& Editing , Visualization. \textbf{Matthias Cosnefroy}: Methodology, Formal analysis, Writing - Review \& Editing. \textbf{Allan Peter Engsig-Karup}: Methodology, Formal analysis, Writing - Review \& Editing, Supervision, Funding acquisition. \textbf{Tim Warburton}: Formal analysis, Writing - Review \& Editing. \textbf{Finnur Pind}: Resources, Funding acquisition. \textbf{Cheol-Ho Jeong}: Methodology, Formal analysis, Writing - Review \& Editing, Supervision, Project administration, Funding acquisition.

\bibliographystyle{elsarticle-num} 
\bibliography{bibliography}

\begin{thebibliography}{10}
\expandafter\ifx\csname url\endcsname\relax
  \def\url#1{\texttt{#1}}\fi
\expandafter\ifx\csname urlprefix\endcsname\relax\def\urlprefix{URL }\fi
\expandafter\ifx\csname href\endcsname\relax
  \def\href#1#2{#2} \def\path#1{#1}\fi

\bibitem{Bayliss1980}
A.~{B}ayliss, E.~{T}urkel, {R}adiation {B}oundary {C}ondition for {W}ave-{L}ike {E}quations, {C}ommunications on {P}ure and {A}pplied {M}athematics 33 (1980) 707 -- 725.

\bibitem{Engquist1979}
B.~{E}ngquist, A.~{M}ajda, {R}adiation {B}oundary {C}onditions for {A}coustic and {E}lastic {W}ave {C}alculations, {C}ommunications on {P}ure and {A}pplied {M}athematics 32 (1979) 313 -- 357.

\bibitem{Collino1993}
F.~{C}ollino, {H}igh {O}rder {A}bsorbing {B}oundary {C}onditions for {W}ave {P}ropagation {M}odels. {S}traight {L}ine {B}oundary and {C}orner {C}ases, {P}roc. 2nd {I}nt. {C}onf. on {M}athematical \& {N}umerical {A}spects of {W}ave {P}ropagation (01 1993).

\bibitem{Appelo2009}
D.~{A}ppelö, T.~{C}olonius, {A} {H}igh-{O}rder {S}uper-{G}rid-scale {A}bsorbing {L}ayer and its {A}pplication to {L}inear {H}yperbolic {S}ystems, {J}ournal of {C}omputational {P}hysics 228 (2009) 4200--4217.

\bibitem{Israeli1981}
M.~{I}sraeli, S.~A. {O}rszag, {A}pproximation of radiation {B}oundary {C}onditions, {J}ournal of {C}omputational {P}hysics 41~(1) (1981) 115--135.

\bibitem{Shlomo1995}
T.~{S}hlomo, N.~{D}ouglas, {A}n {A}bsorbing {B}uffer {Z}one {T}echnique for {A}coustic {W}ave {P}ropagation, 33rd {A}erospace {S}ciences {M}eeting and {E}xhibit (01 1995).
\newblock \href {https://doi.org/10.2514/6.1995-164} {\path{doi:10.2514/6.1995-164}}.

\bibitem{Berenger1994}
J.-P. {B}\'{e}renger, {A} {P}erfectly {M}atched {L}ayer for the absorption of electromagnetic waves, {J}ournal of {C}omputational {P}hysics 114~(2) (1994) 185--200.
\newblock \href {https://doi.org/10.1006/jcph.1994.1159} {\path{doi:10.1006/jcph.1994.1159}}.

\bibitem{Chew1994}
W.~C. {C}hew, W.~H. {W}eedon, {A} 3{D} {P}erfectly {M}atched medium from modified {M}axwell's equations with stretched coordinates, {M}icrowave and {O}ptical {T}echnology {L}etters 7 (1994) 599--604.

\bibitem{Abarbanel1997}
S.~{A}barbanel, D.~{G}ottlieb, {A} {M}athematical {A}nalysis of the {P}{M}{L} {M}ethod, {J}ournal of {C}omputational {P}hysics 134~(2) (1997) 357--363.
\newblock \href {https://doi.org/https://doi.org/10.1006/jcph.1997.5717} {\path{doi:https://doi.org/10.1006/jcph.1997.5717}}.

\bibitem{Hesthaven1998}
J.~S. {H}esthaven, {O}n the {A}nalysis and {C}onstruction of {P}erfectly {M}atched {L}ayers for the {L}inearized {E}uler {E}quations, {J}ournal of {C}omputational {P}hysics 142~(1) (1998) 129--147.
\newblock \href {https://doi.org/https://doi.org/10.1006/jcph.1998.5938} {\path{doi:https://doi.org/10.1006/jcph.1998.5938}}.

\bibitem{Becache2003}
E.~{B}\'{e}cache, S.~{F}auqueux, P.~{J}oly, {S}tability of {P}erfectly {M}atched {L}ayers, group velocities and anisotropic waves, {J}ournal of {C}omputational {P}hysics 188~(2) (2003) 399--433.
\newblock \href {https://doi.org/https://doi.org/10.1016/{S}0021-9991(03)00184-0} {\path{doi:https://doi.org/10.1016/{S}0021-9991(03)00184-0}}.

\bibitem{Diaz2006}
J.~{D}iaz, P.~{J}oly, {A} time domain analysis of {PML} models in acoustics, {C}omputer {M}ethods in {A}pplied {M}echanics and {E}ngineering 195~(29) (2006) 3820--3853.
\newblock \href {https://doi.org/https://doi.org/10.1016/j.cma.2005.02.031} {\path{doi:https://doi.org/10.1016/j.cma.2005.02.031}}.

\bibitem{Hu2001}
F.~{H}u, {A} {S}table, {P}erfectly {M}atched {L}ayer for {L}inearized {E}uler {E}quations in {U}nsplit {P}hysical {V}ariables, {J}ournal of {{C}}omputational {{P}}hysics 173 (2001) 455--480.
\newblock \href {https://doi.org/10.1006/jcph.2001.6887} {\path{doi:10.1006/jcph.2001.6887}}.

\bibitem{Hu2005}
F.~{H}u, {A} {P}erfectly {M}atched {L}ayer absorbing {B}oundary {C}ondition for {LEE} with a non-uniform {M}ean {F}low, {J}ournal of {C}omputational {P}hysics 208 (2005) 469--492.
\newblock \href {https://doi.org/10.1016/j.jcp.2005.02.028} {\path{doi:10.1016/j.jcp.2005.02.028}}.

\bibitem{Parrish2009}
S.~{P}arrish, F.~{H}u, {PML} {A}bsorbing {B}oundary {C}onditions for the {L}inearized and {N}onlinear {E}uler {E}quations in the case of {O}blique {M}ean {F}low, {I}nternational {J}ournal for {N}umerical {M}ethods in {F}luids 60 (2009) 565 -- 589.
\newblock \href {https://doi.org/10.1002/fld.1905} {\path{doi:10.1002/fld.1905}}.

\bibitem{Meza2008}
K.~{M}eza {F}ajardo, A.~{P}apageorgiou, {A} {N}onconvolutional, {S}plit-{F}ield, {P}erfectly {M}atched {L}ayer for {W}ave {P}ropagation in {I}sotropic and {A}nisotropic {E}lastic {M}edia: {S}tability {A}nalysis, {B}ulletin of the {S}eismological {S}ociety of {A}merica 98 (2008) 1811--1836.
\newblock \href {https://doi.org/10.1785/0120070223} {\path{doi:10.1785/0120070223}}.

\bibitem{Roden2000}
J.~{R}oden, S.~{G}edney, {C}onvolution {PML} ({CPML}): {A}n {E}fficient {FDTD} implementation of the {CFS}-{PML} for arbitrary media, {M}icrowave and {O}ptical {T}echnology {L}etters 27 (2000) 334--339.

\bibitem{Grote2010}
M.~J. {G}rote, I.~{S}im, {E}fficient {PML} for the wave equation, ar{X}iv preprint ar{X}iv:1001.0319. (2010).

\bibitem{Etienne2010}
V.~{E}tienne, E.~{C}haljub, J.~{V}irieux, N.~{G}linsky, {A}n hp-adaptive discontinuous {G}alerkin {F}inite-{E}lement {M}ethod for 3-{D} {E}lastic {W}ave {M}odelling, {G}eophysical {J}ournal {I}nternational 183~(2) (2010) 941 – 962.
\newblock \href {https://doi.org/10.1111/j.1365-246{X}.2010.04764.x} {\path{doi:10.1111/j.1365-246{X}.2010.04764.x}}.

\bibitem{Komatitsch2007}
D.~{K}omatitsch, R.~{M}artin, {A}n {U}nsplit {C}onvolutional {P}erfectly {M}atched {L}ayer improved at {G}razing {I}ncidence for the {S}eismic {W}ave {E}quation, Geophysics 72 (09 2007).
\newblock \href {https://doi.org/10.1190/1.2757586} {\path{doi:10.1190/1.2757586}}.

\bibitem{Martin2009}
R.~{M}artin, D.~{K}omatitsch, {{A}n {U}nsplit {C}onvolutional {P}erfectly {M}atched {L}ayer {T}echnique {I}mproved at {G}razing {I}ncidence for the {V}iscoelastic {W}ave 7quation}, Geophysical Journal International 179~(1) (2009) 333--344.

\bibitem{Cosnefroy2019}
M.~Cosnefroy, {S}imulation {N}um{\'e}rique de la {P}ropagation dans l'{A}tmosph{\`e}re de {S}ons {I}mpulsionnels et {C}onfrontations {E}xp{\'e}rimentales [numerical simulation of atmospheric propagation of impulse sound and experimental comparisons], Theses, {E}cole {C}entrale de {L}yon (May 2019).

\bibitem{Becache2021}
E.~{B}\'{e}cache, M.~{K}achanovska, {S}tability and {C}onvergence {A}nalysis of {T}ime-{D}omain {P}erfectly {M}atched {L}ayers for the {W}ave {E}quation in {W}aveguides, {S}{I}{A}{M} {J}ournal on {N}umerical {A}nalysis 59~(4) (2021) 2004--2039.
\newblock \href {https://doi.org/10.1137/20{M}1330543} {\path{doi:10.1137/20{M}1330543}}.

\bibitem{Baffet2019}
D.~{B}affet, M.~{G}rote, S.~{I}mperiale, M.~{K}achanovska, {E}nergy {D}ecay and {S}tability of a {P}erfectly {M}atched {L}ayer {F}or the {W}ave {E}quation, {J}ournal of {S}cientific {C}omputing 81 (12 2019).
\newblock \href {https://doi.org/10.1007/s10915-019-01089-9} {\path{doi:10.1007/s10915-019-01089-9}}.

\bibitem{Kaltenbacher2013}
B.~{K}altenbacher, M.~{K}altenbacher, I.~{S}im, {A} modified and stable version of a {P}erfectly {M}atched {L}ayer technique for the 3-{D} second order wave equation in time domain with an application to aeroacoustics, {J}ournal of {C}omputational {P}hysics 235 (2013) 407--422.
\newblock \href {https://doi.org/10.1016/j.jcp.2012.10.016} {\path{doi:10.1016/j.jcp.2012.10.016}}.

\bibitem{Hesthaven2007}
J.~S. {H}esthaven, T.~{W}arburton, {N}odal {D}iscontinuous {G}alerkin {M}ethods: {A}lgorithms, {A}nalysis, and {A}pplications, 1st Edition, Vol.~54, {S}pringer {P}ublishing {C}ompany, {I}ncorporated, 2007.

\bibitem{ReedHill1972}
W.~H. {R}eed, T.~R. {H}ill, {T}riangular mesh methods for the neutron transport equation, Tech. Rep. ({L}{A}-{U}{R}-73-479), {L}os {A}lamos {S}cientific {L}ab., {N}{M} (October 1973).

\bibitem{Atkins1998}
H.~L. {A}tkins, C.-W. {S}hu, {Q}uadrature-{F}ree {I}mplementation of {D}iscontinuous {G}alerkin {M}ethod for {H}yperbolic {E}quations, {A}{I}{A}{A} {J}ournal 36~(5) (1998) 775--782.
\newblock \href {https://doi.org/10.2514/2.436} {\path{doi:10.2514/2.436}}.

\bibitem{Toulopoulos2006}
I.~{T}oulopoulos, J.~{E}katerinaris, {H}igh-{O}rder {D}iscontinuous {G}alerkin {D}iscretizations for {C}omputational {A}eroacoustics in {C}omplex {D}omains, {A}iaa {J}ournal - {A}{I}{A}{A} {J} 44 (2006) 502--511.
\newblock \href {https://doi.org/10.2514/1.11422} {\path{doi:10.2514/1.11422}}.

\bibitem{Kreiss1972}
H.-O. Kreiss, J.~Oliger, Comparison of accurate methods for the integration of hyperbolic equations, Tellus 24~(3) (1972) 199--215.
\newblock \href {https://doi.org/10.3402/tellusa.v24i3.10634} {\path{doi:10.3402/tellusa.v24i3.10634}}.

\bibitem{Langtangen1998}
H.~P. {L}angtangen, G.~{P}edersen, {C}omputational models for weakly dispersive nonlinear water waves, {C}omputer {M}ethods in {A}pplied {M}echanics and {E}ngineering 160~(3) (1998) 337--358.
\newblock \href {https://doi.org/10.1016/{S}0045-7825(98)00293-{X}} {\path{doi:10.1016/{S}0045-7825(98)00293-{X}}}.

\bibitem{Melander2023}
A.~{M}elander, E.~{S}trøm, F.~{P}ind, A.~P. {E}ngsig {K}arup, C.-H. {J}eong, T.~{W}arburton, N.~{C}halmers, J.~S. {H}esthaven, {M}assively parallel nodal discontinous {G}alerkin finite element method simulator for room acoustics, {T}he {I}nternational {J}ournal of {H}igh {P}erformance {C}omputing {A}pplications (2023) 10943420231208948\href {https://doi.org/10.1177/10943420231208948} {\path{doi:10.1177/10943420231208948}}.

\bibitem{Wang2019}
H.~{W}ang, I.~{S}ihar, R.~P. {M}u{\~n}oz, M.~{H}ornikx, {R}oom acoustics modelling in the time-domain with the nodal discontinuous {G}alerkin method., {T}he {J}ournal of the {A}coustical {S}ociety of {A}merica 145 (2019) 2650.

\bibitem{Wang2020}
H.~{W}ang, M.~{H}ornikx, {{T}ime-{D}omain impedance {B}oundary {C}ondition modeling with the discontinuous {G}alerkin method for room acoustics simulations}, {T}he {J}ournal of the {A}coustical {S}ociety of {A}merica 147~(4) (2020) 2534--2546.
\newblock \href {https://doi.org/10.1121/10.0001128} {\path{doi:10.1121/10.0001128}}.

\bibitem{Pind2020}
F.~{P}ind, C.-H. {J}eong, A.~P. {E}ngsig {K}arup, J.~S. {H}esthaven, J.~{S}trømann {A}ndersen, {T}ime-domain room acoustic simulations with extended-reacting porous absorbers using the discontinuous {G}alerkin method, {T}he {J}ournal of the {A}coustical {S}ociety of {A}merica 148 (2020) 2851--2863.
\newblock \href {https://doi.org/10.1121/10.0002448} {\path{doi:10.1121/10.0002448}}.

\bibitem{Pind2021}
F.~{P}ind, C.-H. {J}eong, J.~S. {H}esthaven, A.~P. {E}ngsig {K}arup, J.~{S}trømann {A}ndersen, {A} phenomenological extended-reaction boundary model for time-domain wave-based acoustic simulations under sparse reflection conditions using a wave splitting method, {A}pplied {A}coustics 172 (2021) 107596.
\newblock \href {https://doi.org/https://doi.org/10.1016/j.apacoust.2020.107596} {\path{doi:https://doi.org/10.1016/j.apacoust.2020.107596}}.

\bibitem{Modave2017}
A.~Modave, J.~{L}ambrechts, C.~{G}euzaine, {P}erfectly {M}atched {L}ayers for {C}onvex {T}runcated {D}omains with {D}iscontinuous {G}alerkin {T}ime {D}omain {S}imulations, Computers and {M}athematics with {A}pplications 73 (01 2017).
\newblock \href {https://doi.org/10.1016/j.camwa.2016.12.027} {\path{doi:10.1016/j.camwa.2016.12.027}}.

\bibitem{LeVeque2002}
R.~J. {L}e{V}eque, {F}inite {V}olume {M}ethods for {H}yperbolic {P}roblems, {C}ambridge {T}exts in {A}pplied {M}athematics, {C}ambridge {U}niversity {P}ress, 2002.

\bibitem{Carpenter1994}
M.~H. {C}arpenter, C.~A. {K}ennedy, {F}ourth-order 2{N}-storage {R}unge-{K}utta schemes, {T}echnical {M}emorandum ({NASA}-{TM}-109112), {NASA} {L}angley {R}esearch {C}enter, {H}ampton, {VA} (June 1994).

\bibitem{Warburton2008}
T.~Warburton, T.~Hagstrom, Taming the cfl number for discontinuous galerkin methods on structured meshes, SIAM J. Numerical Analysis 46 (2008) 3151--3180.
\newblock \href {https://doi.org/10.1137/060672601} {\path{doi:10.1137/060672601}}.

\bibitem{Engsig2006}
A.~P. {E}ngsig {K}arup, {{U}nstructured {N}odal ({DG}-{FEM}) solution of high-order {B}oussinesq-type equations}, Ph.D. thesis, {{T}echnical {U}niversity of {D}enmark ({DTU})} ({A}ugust 2006).

\bibitem{Geuzaine2009Gmsh}
C.~{G}euzaine, J.-F. {R}emacle, {G}msh: {A} 3-{{D}} {F}inite {E}lement {M}esh {G}enerator with built-in {P}re- and {P}ost-{P}rocessing {F}acilities, {I}nternational {J}ournal for {N}umerical {M}ethods in {E}ngineering 79~(11) (2009) 1309--1331.
\newblock \href {https://doi.org/10.1002/nme.2579} {\path{doi:10.1002/nme.2579}}.

\bibitem{Modave2014}
A.~{M}odave, E.~{D}elhez, C.~{G}euzaine, {O}ptimizing {P}erfectly {M}atched {L}ayers in discrete contexts, {I}nternational {J}ournal for {N}umerical {M}ethods in {E}ngineering 99 (08 2014).
\newblock \href {https://doi.org/10.1002/nme.4690} {\path{doi:10.1002/nme.4690}}.

\bibitem{libParanumal}
N.~Chalmers, A.~{K}arakus, A.~P. {A}ustin, K.~{S}wirydowicz, T.~{W}arburton, \href{https://github.com/paranumal/libparanumal}{{libParanumal}: a performance portable high-order finite element library}, release 0.5.0 (2022).
\newblock \href {https://doi.org/10.5281/zenodo.4004744} {\path{doi:10.5281/zenodo.4004744}}.
\newline\urlprefix\url{https://github.com/paranumal/libparanumal}

\end{thebibliography}
\end{document}